\documentclass{amsproc}

\usepackage{graphicx}

\usepackage{hyperref}
\usepackage{latexsym,amsmath,amsfonts,amscd,amssymb, mathdots,bm}
\usepackage{lscape}
\usepackage{array}
\usepackage[all]{xy}
\usepackage{epsfig,color}
\usepackage[version=0.96]{pgf} 

\newtheorem{theorem}{Theorem}[section]

\theoremstyle{definition}
\newtheorem{definition}[theorem]{Definition}
\newtheorem{example}[theorem]{Example}

\newtheorem{corollary}[theorem]{Corollary}

\theoremstyle{remark}
\newtheorem{remark}[theorem]{Remark}

\numberwithin{equation}{section}

\newcommand{\R}{\mathbb{R}}
\newcommand{\Z}{\mathbb{Z}}
\newcommand{\C}{\mathbb{C}}
\newcommand{\noi}{\noindent}
\newcommand{\Mhiggs}{\mathcal{M}(\Sigma,}
\newcommand{\Ksigma}{K_{\Sigma}}
\newcommand{\PU}{\mathrm{PU}}

\newcommand{\PGL}{\mathrm{PGL}}
\newcommand{\PSL}{\mathrm{PSL}}
\newcommand{\PSO}{\mathrm{PSO}}
\newcommand{\SU}{\mathrm{SU}}

\newcommand{\GL}{\mathrm{GL}}
\newcommand{\SL}{\mathrm{SL}}

\newcommand{\SO}{\mathrm{SO}}
\newcommand{\Sp}{\mathrm{Sp}}
\newcommand{\PSp}{\mathrm{PSp}}
\newcommand{\cM}{\mathcal{M}}
\newcommand{\Hom}{\mathrm{Hom}}

\begin{document}

\title{Global properties of Higgs bundle moduli spaces}


\author{Steven Bradlow}
\address{University of Illinois at Urbana-Champaign}
\curraddr{}
\email{bradlow@illinois.edu}
\thanks{}


Key words = Higgs bundles, moduli spaces, connected components

\date{}

\begin{abstract} The moduli spaces for Higgs bundles associated to real Lie groups and a closed Riemann surface have multiple connected components. This survey provides a compendium of results concerning the counting of these components in cases where the  Lie group is a real forms of a complex simple Lie group. In some cases the components can be described quite explicitly. 
\end{abstract}

\maketitle
\tableofcontents

\hspace{.2in}\parbox{4in}{\begin{center}{\it Dedicated with gratitude and admiration to Peter Newstead on the occasion of his 80th birthday}
\end{center}}

\section{Introduction}

The goals of this survey are quite modest: to collect in one place results on the global features of the moduli spaces of Higgs bundles on a Riemann surface, in particular 
\begin{itemize}
\item results on the number of connected components, and
\item descriptions of global structure 
\end{itemize}

The moduli spaces in question, denoted here by $\Mhiggs G)$, are associated to a Riemann surface $\Sigma$ and a Lie group $G$.  We consider only closed surfaces, i.e.\\ we do not allow the surface to have punctures or non-empty boundary.  Many of the basic results referenced in this survey apply if $G$ is any real reductive Lie group\footnotemark\footnotetext{As defined in \cite{Kn}}, but we will assume that $G$ is a non-compact real form of a connected complex semisimple Lie Group.  The Lie algebras of such real forms thus encompass all the non-compact real forms of the complex simple Lie algebras.  

One of the many points of interest concerning the moduli spaces  $\Mhiggs G)$ is the so-called non-Abelian Hodge (NAH) correspondence.  This correspondence identifies $\Mhiggs G)$ with a moduli space that depends only on the underlying real surface of $\Sigma$, say $S$. Denoted by $\mathrm{Rep}(\pi_1(S),G)$, this moduli space parameterizes conjugacy classes of homomorphisms from the fundamental group of $S$ into $G$. Individually, the two moduli spaces have been viewed from very different perspectives,  with each perspective revealing characteristics not easily seen from the other. 

On the Higgs bundle side the moduli space depends on choices, specifically a choice of complex structure on the surface $S$. The reward for this arbitrary choice is that the Higgs moduli space is endowed with extra structures which serve as tools for exploring the moduli space. Perhaps surprisingly, these tools reveal properties which are independent of the complex structure and hence transfer to the representation variety side of the correspondence. 

The correspondence relies on existence and uniqueness results for non-linear PDE's. As a result it is very difficult to pair individual Higgs bundles with specific representations or vice versa.    The information gleaned on the representation varieties side of the NAH correspondence leans towards the properties of the individual representations in specific components, or the dynamics of the mapping class group action, with particularly striking results for so-called higher Teichm\"uller  components. In contrast, the Higgs bundle tools are very effective at revealing topological information such as the number of connected components and, in some cases, other global features.  

We will remain firmly on the Higgs bundle side of the NAH correspondence.  There are by now several excellent surveys of $G$-Higgs bundles on Riemann surfaces and their moduli spaces (see Section \ref{Ghiggsbundles}). We will lean heavily on these and provide here a much abbreviated summary of the main points relevant for our purposes. 

In keeping with our limited goals there are many key topics that are not covered in this survey.  In particular we will say hardly anything about the many interesting aspects of the representation varieties $\mathrm{Rep}(\pi_1(S), G)$ and their components. On the Higgs bundle side of the NAH correspondence we will also be forced to give short shrift to many important topics, including 
\begin{itemize}
\item Spectral data for $G$-Higgs bundles when $G$ is a real form, and the use of monodromy in the Hitchin fibration to count components. The Hitchin fibration as described in Section \ref{HitchinFibComponents} has singular fibers. Removing these leaves a true non-singular fibration over the regular locus in the base of the fibration. The monodromy around the singular fibers, in conjunction with spectral data for real Higgs bundles can be used to infer an upper bound on the number of connected components in $\Mhiggs G)$.  If this agrees with a known lower bound then the method can provide an elegant way to count components.  While spectral data for real Higgs bundles add significant extra facets to our understanding (see for example \cite{BS19}), to date this method has provided new information about the number of connected components only in the cases of $\mathrm{GL}(2,\R)$ and $\SO(2,2)$ (see \cite{BSrank2}).
\item Geometric structures detected by $G$-Higgs bundles. One of the many interesting features of so-called higher Teichm\"uller  components of the representation varieties $\mathrm{Rep}(\pi_1(S),G)$ is their relation to geometric structures either on $S$ itself or on related spaces (see \cite{W2018} or most recently \cite{AMTW}). This may be viewed as a generalization of the way the Teichm\"uller  space of $S$ appears as a component of $\mathrm{Rep}(\pi_1(S),PSL(2,\R))$. Such generalizations have been explored mostly from the point of view of the representations, but $G$-Higgs bundles have increasingly been used in interesting new ways (see \cite{allesandrini2019} for a recent survey).
\end{itemize}

Many of the results covered by this survey resulted initially from a case-by-case analysis, starting with Hitchin's groundbreaking study of $\SL(2,\R)$-Higgs bundles in \cite{selfduality}. Indeed the techniques introduced there and in \cite{liegroupsteichmuller} permeate virtually all of the subsequent developments. This is certainly true for the many specific examples of non-compact real forms that have been individually studied (see the references in the Tables), but it is also true for the subsequent general frameworks that have followed. In addition to Hitchin's results for split real forms in \cite{liegroupsteichmuller}, these include a general result for real forms of Hermitian type \cite{BGRmaximalToledo} and, most recently, the construction in \cite{BCGGO21} of so-called Cayley components of $\Mhiggs G)$ for real forms arising from distinguished (`magical') $\mathfrak{sl}(2,\C)$-subalgebras.

Rather than starting with the most general results and then working out some noteworthy special cases, we will proceed in the opposite direction. The goal is to build intuition to make it easier to appreciate the mechanisms at work in the general unifying result.  We hope to strike a balance between providing all the details and getting to the endpoint in a timely manner.

\medskip
\noi {\bf Acknowledgement.} There are many colleagues on both sides of the non-Abelian Hodge Correspondence, too many to list,  who have immeasurably enriched my mathematical life in general and in particular my understanding of the work discussed here. Foremost among them are my main collaborators Oscar Garcia-Prada, Brian Collier, Peter Gothen, and Andre Oliveira,  but others to whom I owe a great debt of gratitude include Bill Goldman, Nigel Hitchin, and Anna Wienhard. My sincere thanks to the editors of this volume for their truly Job-level of patience, and to the anonymous referee whose sharp insights forced clarifications of several poorly explained points and whose extraordinary diligence sanded down numerous rough edges in the text.

 \section{The real forms}
 
Real forms of a complex simple Lie group $G^{\C}$ are by definition the fixed points of an anti-holomorphic Lie group involution on $G^{\C}$. These always include a compact real form, unique up to conjugation, for which we will denote the anti-holomorphic involution by
\begin{equation}
\kappa: G^{\C}\rightarrow G^{\C}\ .
\end{equation}
If $\tau: G^{\C}\rightarrow G^{\C}$ is any other anti-holomorphic involution, then $\sigma=\tau\kappa$ defines a holomorphic involution.  The holomorphic involution is called the {\it Cartan involution} for $\tau$.  If the compact real form is fixed, i.e.\ if $\kappa$ is given, then the non-compact real forms are specified either by the other anti-holomorphic involutions or their holomorphic Cartan involutions.

The holomorphic and anti-holomorphic involutions on the Lie groups induce corresponding linear or conjugate-linear involutions on their Lie algebras. The real forms of the complex Lie algebra are the fixed point sets of the conjugate-linear involutions.

The real groups of interest in this survey are non-compact real forms of complex simple Lie groups. As such their Lie algebras are non-compact real forms of complex simple Lie algebras. In Table \ref{realforms} we list all such real forms organized according to features which play a role in the rest of this survey, i.e.\ which are important for understanding the components of the Higgs bundle moduli spaces.  

Some properties of the Higgs bundle moduli spaces depend only on Lie algebra data but others, for instance in the precise count of their connected components, are sensitive to the center of the group  - a feature not detected by the Lie algebra.  Among the connected Lie groups with a given Lie algebra there is always an adjoint form, i.e.\ a Lie group with trivial center.  The adjoint form has a universal cover which is also a Lie group with the same Lie algebra.  All the other groups are quotients of the universal cover and are finite covers of the adjoint group.  Note that not all such covers are matrix groups. For example the finite covers of $\SL(2,\R)$ cannot be realized as matrix groups. Though Higgs bundles can in principle be defined for all such finite covers we will consider only matrix groups. 

Another feature of a group that is not determined by its Lie algebra is its number of connected components. Some real forms, notably $\SO(p,q)$, have more than one component. For these real forms the component of the identity (e.g. $\SO_0(p,q)$ in the case of $\SO(p,q)$) is still a real reductive Lie group as defined in \cite{Kn}. Our survey includes such cases. 

\clearpage
\pagebreak

\begin{table}
\centering
  \resizebox{\textwidth}{!}{  
\begin{tabular}{|c||c|c|c||c|c|}
\hline
$\mathfrak{g}^{\C}$ & Split & Hermitian & Other& Hermitian  & Other\\
 & &tube &Magical& non-tube &non-Magical\\
\hline
 & & & & &\\
$\mathfrak{sl}(n,\C)$ &$\mathfrak{sl}(n,\R)$ & $\mathfrak{su}(p,p)$ & --&$\mathfrak{su}(p,n-p)$ & $\mathfrak{su}^*(2p)$\\
& & ($n=2p$) &&$1\le p<\frac{n}{2}$  &($n=2p$)\\
  \hline
 & & & & &\\
  $\mathfrak{so}(2n+1,\C)$ &$\mathfrak{so}(n,n+1)$ & $\mathfrak{so}(2,2n-1)$ & $\mathfrak{so}(p,2n+1-p)$ &--& $\mathfrak{so}(1,2n)$\\
  & & & $2<p\le n$ &&\\
 \hline
 & & & & &\\
 $\mathfrak{sp}(2n,\C)$ &$\mathfrak{sp}(2n,\R)$ &$\mathfrak{sp}(2n,\R)$ &--&--& $\mathfrak{sp}(2p,2n-2p)$\\
  & & & & &($1\le p\le \frac{n}{2}$)\\
  \hline
 & & & & &\\
  $\mathfrak{so}(2n,\C)$ &$\mathfrak{so}(n,n)$ & $\mathfrak{so}(2,2n-2)$ & $\mathfrak{so}(p,2n-p)$ &$\mathfrak{so}^*(4p+2)$&  $\mathfrak{so}(1,2n-1)$\\
  & &$\mathfrak{so}^*(4p) (n=2p)$ &$2<p < n$ &$(n=2p+1)$ &\\
 \hline
 & & & & &\\
 $\mathfrak{g}_2$  & $\mathfrak{g}_2^2$& --&-- &-- & --\\
 & & & & &\\
 \hline
 & & & & &\\
 $\mathfrak{f}_4$  &$\mathfrak{f}_4^4$  & --&  $\mathfrak{f}_4^4$ &--& $\mathfrak{f}_4^{-20}$ \\
 & & & & &\\
  \hline
 & & & & &\\
 $\mathfrak{e}_6$  &$\mathfrak{e}_6^6$ &-- & $\mathfrak{e}_6^{2}$&$\mathfrak{e}_6^{-14}$ &$\mathfrak{e}_6^{-26}$\\
 & & & & &\\
  \hline
 & & & & &\\
 $\mathfrak{e}_7$ & $\mathfrak{e}^7_7$&$\mathfrak{e}^{-25}_7$ & $\mathfrak{e}^{-5}_7$&--& --\\
 & & & & &\\
  \hline
 & & & & &\\
 $\mathfrak{e}_8$  & $\mathfrak{e}^8_8$& --&$\mathfrak{e}^{-24}_8$ &--&--\\
 & & & & &\\
   \hline
\end{tabular}
}
 \medskip
\caption {Non-compact real forms of complex simple Lie algebras (with $\mathfrak{sp}(2n,\R), \mathfrak{so}(n,n+1),\mathfrak{f}_4$ listed twice because they have two magical structures)}\label{realforms}
\end{table}

\clearpage
\pagebreak

\section{$G$-Higgs bundles}\label{Ghiggsbundles}

There are several recent surveys of Higgs bundles on Riemann surfaces and their moduli spaces (see \cite{garcia-pradasurvey09, garcia-pradasurvey, gothensurvey, guichardsurvey, Li2019, schaposniksurvey, swobodasurvey, wentworthsurvey}, and also Appendix A in \cite{so(pq)BCGGO}) so we will repeat here only the bare essentials. 

Let $S$ be a closed oriented surface of genus $g\ge 2$. We fix a complex structure on $S$ and denote the resulting Riemann surface by $\Sigma$. We also fix a compatible Riemannian metric on $\Sigma$. In general terms, a $G$-Higgs bundle on $\Sigma$ is a pair $(E,\Phi)$ where $E$ is a holomorphic bundle on $\Sigma$ and $\Phi$, the Higgs field, is a holomorphic section of an associated vector vector bundle twisted by the canonical bundle on $\Sigma$.  If $G$ is a complex Lie group, then $E$ is a holomorphic principal $G$-bundle and the associated vector bundle is the adjoint bundle $ad(E)$, i.e.\ a bundle with Lie algebra of $G$ as fiber. The Higgs field is then a section of $ad(E)\otimes \Ksigma$, where $\Ksigma$ is the canonical bundle.  For our present purposes we need a more flexible framework in which $G$ is a real Lie group. In this setting the description of both $E$ and the bundle containing the Higgs field is less immediate than in the complex case.  In all cases, the definition is motivated by properties described in the next sections.

\subsection{Definition}\label{GHiggsDefn}

Let  $G$ be a non-compact real form of a complexl semisimple\footnotemark\footnotetext{the framework may straightforwardly be enlarged to include reductive Lie groups but we will confine attention to the semisimple case} Lie group, $G^{\C}$. Let $H\subset G$ be a maximal compact subgroup of $G$ and denote by $H^{\C}$ its complexification.  The Cartan involution defining $G$ induces a decomposition of the Lie algebra $\mathfrak{g}^{\C}$ of $G^{\C}$

\begin{equation}\label{gcdecomp}
\mathfrak{g}^{\C}=\mathfrak{h}^{\C}\oplus \mathfrak{m}^{\C}\ ,
\end{equation}

\noi where the summands are the $\pm 1$-eigenspaces. The $+1$-eigenspace $\mathfrak{h}^{\C}$ is the Lie algebra of $H^{\C}\subset G^{\C}$ where $H^{\C}$ is the complexification of a maximal compact subgroup $H\subset G$.
The adjoint action of $H^{\C}$ preserves $\mathfrak{m}^{\C}$ and the resulting representation of $H^{\C}$ on $\mathfrak{m}^{\C}$ is the so-called {\it isotropy representation}

\begin{definition}\label{GHiggs}A $G$-Higgs bundle on a Riemann surface $\Sigma$ is a pair $(E_{H^{\C}},\phi)$ where $E_{H^{\C}}$ is a holomorphic principal $H^{\C}$-bundle and $\Phi$ is a holomorphic section of the vector bundle $E_{H^{\C}}[\mathfrak{m}^{\C}]\otimes \Ksigma$. Here $E_{H^{\C}}[\mathfrak{m}^{\C}]$ denotes the associated vector bundle (via the isotropy representation) with fiber $\mathfrak{m}^{\C}$, and $\Ksigma$ is the canonical bundle.

Two such Higgs bundles, say $(E_{H^{\C}},\phi)$  and $(E'_{H^{\C}},\phi')$, are isomorphic if there is a bundle isomorphism $u:E_{H^{\C}}\rightarrow E'_{H^{\C}}$ which pulls back $\Phi'$ to $\Phi$, i.e.\ such that $u^*\Phi'=\Phi$.

\end{definition}

The holomorphic bundle $E_{H^{\C}}$ may be viewed as an underlying smooth bundle, say $\mathbf{E}_{H^{\C}}$, together with a $\overline{\partial}$-operator $\overline{\partial}_E$ (also known as an anti-holomorphic partial connection because they arise as the anti-holomorphic part a connection on $\mathbf{E}_{H^{\C}}$). From this point of view a $G$-Higgs bundle is defined as follows.

\begin{definition}\label{GHiggsdbar}A $G$-Higgs bundle on $\mathbf{E}_{H^{\C}}$ is a pair $(\overline{\partial}_E, \Phi)$ where $\Phi$ is  smooth holomorphic 1-forms with values in $\mathbf{E}_{H^{\C}} [\mathfrak{m}^{\C}]$, i.e.\ $\Phi\in \Omega^{1,0}(\mathbf{E}_{H^{\C}}[ \mathfrak{m}^{\C}])$ with  $\overline{\partial}_E \Phi=0$. \end{definition}

The complex gauge group for $\mathbf{E}_{H^{\C}}$, i.e.\ $\mathcal{G}_{H^{\C}}=\Omega^0(Ad(\mathbf{E}_{H^{\C}}))$, acts on the space of pairs  $(\overline{\partial}_E, \Phi)$ with the action preserving the holomorphicity condition, and such that the gauge orbits correspond to isomorphism classes of Higgs bundles as defined in Definition \ref{GHiggs}).

If $G= \SL(n,\C)$ then the principal bundle can be viewed as a frame bundle and the definition is equivalent to the more familiar definition in terms of vector bundles, i.e.\ a rank $n$ Higgs bundle is a pair $(E,\phi)$ in which $E$ is a rank $n$ holomorphic vector bundle with trivial determinant and the Higgs field is a holomorphic section of $End(E)\otimes \Ksigma$. 

The motivation for the definition comes in part from the following considerations
\begin{itemize}
\item the existence theorem (described in \S \ref{eqtnsetc}) for solutions to natural gauge-theoretic equations yields flat connections with holonomy in $G$, and hence a correspondence between $G$-Higgs bundles and  surface group representations in $G$. 
\item good moduli spaces can be constructed parameterizing isomorphism classes of $G$-Higgs bundles.
\item if $\sigma$ is the involution which defines $G$ as a real form of $G^C$ then $\sigma$ induces an involution on the moduli space of $G^C$-Higgs bundles, with the moduli space of $G$-Higgs bundles in the fixed point locus.
\end{itemize}

\subsection{Equations}\label{eqtnsetc}\hfill

The defining data for a $G$-Higgs bundle, including the Lie theoretic structures described above, facilitates the formulation of a system of PDE's which define a metric, $h$, on the principal bundle $E_{H^{\C}}$. Here a metric means a section of the associated bundle $E_{H^{\C}}(H^{\C}/H)$, where $H\subset H^{\C}$ denotes a maximal compact subgroup.  Such a metric determines a unique connection on $E_{H^{\C}}$  --  called the Chern connection -- defined by compatibility with both the metric and the holomorphic structures.  The metric reduces the structure group of $E_{H^{\C}}$ to $H$ and thereby permits an extension of $\tau$, the antiholomorphic involution which defines $H$, to a globally defined map. Combined with conjugation on forms, this yields a map

$$\tau_h: \Omega^{1,0}(E_{H^{\C}}[\mathfrak{m}^{\C}])\rightarrow \Omega^{0,1}(E_{H^{\C}}[\mathfrak{m}^{\C}])\ .$$

\noi The Higgs bundle equations takes the form

\begin{equation}\label{higgseqtn}
F_h-[\Phi,\tau_h(\Phi)]=0
\end{equation}

\noi where $F_h$ denotes the curvature of the Chern connection on $E_{H^{\C}}$.  In the case where $G=\SL(n,\C)$ so that $h$ is just a Hermitian metric on the holomorphic vector bundle $E$, the curvature term is the curvature of the usual Chern connection, and $\tau_h$ is the map from $\Omega^{1,0}(End(E))\rightarrow \Omega^{0,1}(End(E))$ defined by $\tau_h(\Phi)=-\Phi^{*_{h}}$ (where $*_{h}$ denotes the adjoint with respect to the metric $h$ combined with conjugation on (1,0)-forms).

It is sometimes useful to fix a metric $h$ on the smooth principal bundle $\mathbf{E}_{H^{\C}}$, i.e.\ fix a reduction to an $H$-bundle, say $\mathbf{E}_{H}$, and regard \eqref{higgseqtn} as an equation for a pair $(D_h,\Phi)$ where $D_h$ is a connection on $\mathbf{E}_{H}$, and $\Phi$ is a smooth section of the bundle $\mathbf{E}_{H}[\mathfrak{m}^{\C}]\otimes \Ksigma$.  Equivalently, with $\overline{\partial}_E$ given by the antiholorphic part of $D_h$,  the equation may be viewed as an equation for a pair $(\overline{\partial}_E, \Phi)$ as in Definition \ref{GHiggsdbar}. Equation \eqref{higgseqtn} must then be complemented by the holomorphicity condition 

\begin{equation}\label{holeqtn}
\overline{\partial}_E\Phi=0
\end{equation}

\subsection{Stability and moduli spaces}\label{stability} The concept of stability comes from Geometric Invariant Theory and plays a pivotal role in the construction of moduli spaces. It is straightforward to define stability for Higgs vector bundles but more complicated for general $G$-Higgs bundles.  

\begin{definition} A rank $n$ trivial determinant Higgs vector bundle $(E,\phi)$ on $\Sigma$ is {\it semi-stable} if 

\begin{equation}\label{sstability}
\frac{\deg(E')}{\mathrm{rank}(E')} \le 0
\end{equation}

\noi for all $\phi$-invariant subbundles $E'$, i.e.\ for all subbundles such that $\phi(E')\subset E'\otimes \Ksigma$. It is {\it stable} if the inequality is strict for all non-trivial proper $\phi$-invariant subbundles, and {\it polystable} if it decomposes as a direct sum of stable Higgs bundles of lower rank.
\end{definition}

The definitions in the general setting of $G$-Higgs bundles is complicated by the more elaborate machinery needed to describe subobjects of principal bundles and to formulate the numerical replacement for the degree of a vector bundle.  The details will not be needed in this survey.  A good account can be found in \cite{garcia-pradasurvey}.  The crucial consequence is that the set of isomorphism classes of polystable Higgs bundles can be endowed with the structure of a moduli space. 

\begin{definition} The moduli space of $G$-Higgs bundles on $\Sigma$, denoted $\Mhiggs G)$, is the set of isomorphism classes of polystable $G$-Higgs bundles. 
\end{definition}

\begin{remark} In Definition \ref{GHiggs} if the canonical bundle is replaced by any other line bundle on $\Sigma$, say $L$, we refer to the resulting pairs as {\it $L$-twisted $G$-Higgs pairs}.  The stability concepts extend readily to these more general objects and permit the construction of moduli spaces. We will denote such moduli spaces by $\mathcal{M}_L(\Sigma,G)$. The moduli spaces of  $L$-twisted $G$-Higgs pairs (with $L$ being a power of $\Ksigma$) will play a significant role in Sections \ref{maxcomponents} and \ref{magicalSL2}.
\end{remark}

The moduli spaces as defined above can be endowed with the structure of complex analytic varieties (see for example \cite{schmitt} for an account).  This can also be seen analytically because of the relation between (poly)stability and the existence of a solution to the Higgs bundle equation. 

\begin{theorem}\label{HKforHiggs}(Hitchin-Kobayashi correspondence) The $G$-Higgs bundle $(E_{H^{\C}},\Phi)$ admits a metric satisfying the Higgs bundle equation if and only if $(E_{H^{\C}},\Phi)$ is polystable.
\end{theorem}

Both the existence of solutions to the Higgs bundle equations and also the property of polystability are properties of isomorphism classes, i.e.\ these properties are preserved by $H^{\C}$-gauge transformations. Thus we can idenitfy

\begin{equation}
\Mhiggs G)\simeq \big\{ \big (\overline{\partial}_E, \Phi) \ |  \overline{\partial}_E\Phi=0\ \mathrm{and}\ \eqref{higgseqtn}\ \mathrm{has\ a\ solution}  \}/\mathcal{G}_{H^{\C}}
\end{equation}

\subsection{The non-Abelian Hodge Correspondence}

By extension of structure group $E_{H^{\C}}$ defines a holomorphic $G^{\C}$-bundle which we denote by $E_{H^{\C}}(G^{\C})$ or simply $E_{G^{\C}}$. The Chern connection on $E_{H^{\C}}$ induces a connection on $E_{G^{\C}}$, denoted by $D_h$. Using  \eqref{gcdecomp} we can now define a new connection

\begin{equation}\label{higgscon}
\nabla_h=D_h+\Phi+\Phi^{*_h}
\end{equation}


This is, a priori, a $G^{\C}$-connection, i.e.\ an operator with values in the Lie algebra $\mathfrak{g}^{\C}$.  However if $h$ is a solution to \eqref{higgseqtn} then $\nabla_h$ has two important features:
\begin{itemize}
\item $\nabla_h$ takes values in $\mathfrak{g}$ and
\item $\nabla_h$ is a flat connection, i.e.\ it has zero curvature. 
\end{itemize}

\noi In particular, the holonomy of $\nabla_h$ defines a reduction of structure to $G$ and defines a flat structure on the resulting $G$-bundle. The holonomy of this flat structure defines a representation of $\pi_1(\Sigma)$ in $G$.  In this way, solutions to \eqref{higgseqtn} define representations of $\pi_1(S)$ in $G$.  Conversely, any representation $\rho:\pi_1(S)\rightarrow G$ defines a flat principal $G$-bundle on $S$ via the construction

\begin{equation}
E_G=\tilde{S}\times_{\rho}G\ .
\end{equation}

\noi Here $\tilde{S}$ is the universal cover of $S$, viewed as a principal $\pi_1(S)$-bundle on $S$.

A theorem of Corlette (and also Donaldson)  leads to the following fundamental relation

\begin{theorem}\label{corlette}There is a bijective correspondence between gauge equivalence classes of principal $G$-bundles on $S$ with reductive flat structures and $G$-conjugacy classes of reductive representations in $Hom(\pi_1(S),G)$. 
\end{theorem}

The space of $G$-conjugacy classes of reductive representations in $Hom(\pi_1(S),G)$ is denoted by $\mathrm{Rep}(\pi_1(S),G)$.  The combination of Corlette's theorem and the Hitchin-Kobayashi correspondence for Higgs bundles thus leads to an identification of moduli spaces:

\begin{theorem} Let $\Sigma$ be a closed smooth Riemann surface of genus $g\ge 2$ with underlying topological surface $S$, and let $G$ be a real form of a complex semisimple Lie group. Then 

\begin{equation}\label{NAH}
\Mhiggs G)\simeq\mathrm{Rep}(\pi_1(S),G)
\end{equation}
\end{theorem}

This is known as the non-Abelian Hodge (NAH) correspondence.  We emphasize that the moduli space  $\mathrm{Rep}(\pi_1(S),G)$ depends only on the topology of $S$ (i.e.\ its genus) and on the Lie group $G$, and parameterizes conjugacy classes of surface group representations. In contrast, $\Mhiggs G)$ depends not only on $G$ but on the complex structure on the surface, and parameterizes $G$-Higgs bundles. 

\section{Component types and labels}\label{labels}

On both sides of the NAH correspondence there are invariants coming directly from the topology of $G$.  
If $G$ is connected then in $ \mathrm{Rep}(\pi_1(S),G)$ such invariants arise from considering whether representations into $G$ can be lifted into the universal cover of $G$. These obstruction can be interpreted as an element in $\pi_1(G)$.  In  $\Mhiggs G)$ the topological type of $E_{H^{\C}}$ is constant on connected components.  The topological types of the principal bundles are classified by characteristic classes in $\pi_1(H)=\pi_1(G)$.  If $G$ is not connected the specification of topological types can be more subtle, but in all cases we get a partition of $\Mhiggs G)$ into sectors labelled by the topological types for principal $H^{\C}$-bundles, say
 
 \begin{equation}
 \Mhiggs G)= \coprod_{c\in \mathcal{I}}\mathcal{M}_c(\Sigma, G)
 \end{equation}
  
 \noi where $\mathcal{I}$ denotes the indexing set for the topological types and $\mathcal{M}_c(\Sigma, G)$ denotes the union of connected components in which the principal $H^{\C}$-bundles are of type $c\in \mathcal{I}$. 

\begin{definition}The primary topological invariant of a connected component of $\Mhiggs G)$ is the topological type of the principal bundles in the Higgs bundles parameterized by that component.
\end{definition}

If $G$ is a complex simple group then (see \cite{Oliveira_GarciaPrada_2016}) there is a bijective correspondence between the topological types of $G$-bundles and the connected components of $\Mhiggs G)$. This is not true for the components of $\Mhiggs G)$ if $G$ is a real group. The relation is complicated by factors including
\begin{itemize}
\item (Bounds) Not all topological types may occur. If $G$ is a real form of Hermitian type (see Section \ref{maxcomponents}) then the primary topological invariants determine a discrete real-valued invariant (the Toledo invariant). There is a (Milnor-Wood) bound on the Toledo invariant outside of which the components are empty.

\item (Secondary invariants) The primary invariants may not distinguish connected components. If $G$ admits magical structures (see section \ref{magicalSL2}) then $\Mhiggs G)$ has components in which secondary topological structures emerge. The topological types of these secondary structures are needed to distinguish components with the same primary invariants. 
\end{itemize}

In the complex case, the Higgs bundles can always be deformed to a Higgs bundle in which the Higgs field vanishes. On the other side of the NAH correspondence this corresponds to the fact that the representations can be deformed to representations which factor through the compact real form.  In the next sections we see how the situation differs in the components of the moduli spaces for non-compact real forms, and how this affects the total number of connected components.

The rich structure of $\Mhiggs G)$ offers effective tools  for detecting and counting its connected components. These  include a real-valued function which provides Morse-theoretic information, and direct constructions (of various degrees of explicitness) of special components.  Our main focus in the next sections is to describe such methods and the results they have produced.


 \section{Morse-theoretic approach} \label{morse}

 Let $(E_{H^{\C}},\Phi)$ be a polystable $G$-Higgs bundle. The metric defined by the Higgs bundle equations \eqref{higgseqtn} reduces the structure group of $E_{H^{\C}}$ to a maximal compact subgroup $H\subset H^{\C}$ and hence induces a unitary structure on the bundle $E_{H^{\C}}[\mathfrak{m}^{\C}]=E_H[\mathfrak{m^{\C}}]$. Combined with a fixed metric on the surface, and hence on $\Ksigma$, this permits computation of a fiberwise norm $|\Phi(x)|_h$ at each point $x\in\Sigma$, and hence defines
 
 \begin{equation}\label{fdef}
 f(E_{H^{\C}},\Phi)=||\Phi||_h^2=\int_{\Sigma}|\Phi|_h^2dvol\ .
 \end{equation}
 
 \noi This is constant on isomorphism classes or, equivalently on $H^{\C}$-gauge orbits, and hence produces a well-defined map
 
  \begin{align}\label{fmap}
   f:\Mhiggs G)&\rightarrow\R\\
 [E_{H^{\C}},\Phi]&\mapsto ||\Phi||_h^2\nonumber
 \end{align}

We will call this the Hitchin function as it was first defined by Hitchin who showed:
 
 \begin{theorem}This map is smooth on the smooth locus of the moduli space and is a proper map.
 \end{theorem}
 
 The Hitchin function serves as a Morse-Bott function on the smooth locus in the moduli space. While the singularities in the moduli space complicate this relation, the properness of the function implies that $f$ attains local minima on all connected components of the moduli space. Thus if $f_{min}\subset\Mhiggs G)$ denotes the locus of local minima, then the number of connected components of $\Mhiggs G)$ is bounded above by the number of connected components of $f_{min}$. Moreover the component counts are equal if the connected components of $f_{min}$ all lie in distinct connected components of the moduli space.

Note that Higgs bundles of the form $(E_{H^{\C}},0)$, i.e with $\Phi=0$, necessarily define local minima of $f$. We call these the {\it trivial local minima}.  Crucially for the existence of many connected components, the trivial minima are not necessarily the only local minima. 

One way to identify the non-trivial local minima of the Hitchin function is to use the relation between $f$ and another key feature of the moduli spaces, namely a  $\C^*$-action defined by
 
 \begin{equation}\label{c*action}
 \lambda[E_{H^{\C}},\Phi]=[E_{H^{\C}},\lambda\Phi]
 \end{equation}

 \begin{theorem}The critical points of the Hitchin function coincide with the fixed points of the $\C^*$ action defined by \eqref{c*action}.
 \end{theorem}
 
 \noi We now examine the two types of minima, namely minima with $\Phi=0$ and minima with non-trivial Higgs field.  
 
\subsection{Minima with $\Phi=0$} Suppose $(E_{H^{\C}},0)$ represents a point in  $\Mhiggs G)$. Then $(E_{H^{\C}},0)$ must be polystable as a $G$-Higgs bundle. It follows that $E_{H^{\C}}$ must be polystable as a principal $H^{\C}$ bundle and must support a flat Hermitian Einstein metric, i.e.\ a reduction of structure to $H$ with respect to which the Chern connection has zero curvature. This is not always possible. Depending on the topological type of $E_{H^{\C}}$ there may be obstructions to the existence of flat structures. 

\begin{example} Let $G=\SL(2,\R)$. In this case $H^{\C}=\C^*$ and the topological type of an $\SL(2,\R)$-Higgs bundle is determined by an integer corresponding to the degree of a complex line bundle. All line bundles are stable, but only the zero degree line bundles admit flat structures. 
\end{example}

Components which contain trivial minima are distinguished by the topological types of the principal $H^{\C}$-bundles at the minima.  The $G$-Higgs bundles parameterized by these components can be deformed to $G$-Higgs bundle where the Higgs field is identically zero. On the other side of the NAH correspondence, the representations parameterized by the corresponding components of $\mathrm{Rep}(\pi_1(S),G)$ can be deformed to representations which factor through a maximal compact subgroup of $G$.

However, as the above example illustrates, not all such topological types occur as labels of components with trivial minima. 
 
 \subsection{Minima with non-trivial $\Phi$}   It is important to note that the $\C^*$-action in \eqref{c*action} is on the moduli space, i.e.\ on isomorphism classes of $G$-Higgs bundles.  The condition for an isomorphism class $[E_{H^{\C}},\Phi]\in \Mhiggs,G)$ to be a fixed point is thus that for each $\lambda\in \C^*$ the Higgs bundle $(E_{H^{\C}},\lambda\Phi)$ is isomorphic to $(E_{H^{\C}},\Phi)$. Representing the Higgs bundles by pairs $(\overline{\partial}_E, \Phi)$ as in Definition \ref{GHiggsdbar} this means that  for all $\lambda\in\C^*$ there is a gauge transformation $g(\lambda)\in\Omega^0(E_{H^{\C}}(H^{\C}))$ such that
 
 \begin{equation}\label{c*fixed}
(\overline{\partial}_E, \lambda\Phi)= g(\lambda)^*(\overline{\partial}_E, \Phi),
 \end{equation}
 
\noi  i.e.\ such that
 
 \begin{align}
g(\lambda)^{-1}\circ\overline{\partial}_E\circ g(\lambda)&=\overline{\partial}_E\ \mathrm{and}\label{fixedE}\\
g(\lambda)^{-1}\circ\Phi\circ g(\lambda)&=\lambda\Phi\label{fixedPhi}\ .
 \end{align}

If  the fixed point $[E_{H^{\C}},\lambda\Phi]$ is stable then $ g(\lambda)$ is uniquely determined\footnotemark\footnotetext{For further details see the summary in Appendix A in \cite{so(pq)BCGGO} or for a more algebraic treatment see \cite{BGRmaximalToledo} for a discussion based on \cite{Sim88} of Hodge bundles in the framework of $G$-Higgs bundle.}.  The infinitesimal action at the fixed point induces weight-space decompositions

\begin{align}
E_{H^{\C}}[\mathfrak{h}^{\C}]=\bigoplus_j E_{H^{\C}}[\mathfrak{h}^{\C}]_j\\
E_{H^{\C}}[\mathfrak{m}^{\C}]=\bigoplus_j E_{H^{\C}}[\mathfrak{m}^{\C}]_j
\end{align}

\noi where the weights are integers or half integers.  With respect to these decompositions  $\Phi$ lies in $E_{H^{\C}}[\mathfrak{m}^{\C}]_1\otimes \Ksigma$ and its adjoint action increase weights by one, i.e.\

\begin{equation}\label{adPhi}
ad_{\Phi}: E_{H^{\C}}[\mathfrak{h}^{\C}]_j\rightarrow E_{H^{\C}}[\mathfrak{m}^{\C}]_{j+1}\otimes \Ksigma
 \end{equation}

We are primarily interested in the fixed points corresponding to local minima of the Hitchin function.  Such fixed points may occur within the smooth locus of $\Mhiggs G)$ or at singular points. The latter occur at points represented by strictly polystable Higgs bundles, while the former are all represented by stable Higgs bundles. The following criterion (first introduced in \cite{liegroupsteichmuller}) is the main tool for identifying the smooth local minima.

\begin{theorem}\label{mincrit} Let $[E_{H^{\C}},\Phi]$ be a $\C^*$-fixed point in the smooth locus of $\Mhiggs G)$.  Then $[E_{H^{\C}},\Phi]$ is a local minimum of $f$ if and only if either $\Phi$ is the zero section or the maps in \eqref{adPhi} are sheaf isomorphisms for every $j>0$.
\end{theorem}

The criteria from theorem \ref{mincrit} do not apply to fixed points outside the smooth locus of the moduli space.  In all cases that have been examined so far, non-smooth critical points are analyzed by constructing explicit deformations which decrease the Hitchin function, thereby showing that they are not local minima.

The local minima of the Hitchin function can detect components but they provide no information about the global structure of the components. In the next sections we describe methods which reveal entire components. In some cases this is in the form of explicit descriptions of certain components, while in other cases the information is less complete.

\section{Sections of Hitchin fibrations and Hitchin components}\label{HitchinFibComponents}
 
If  $G$ is a complex semisimple Lie group then the Higgs field in the pair $(E_G,\Phi)$ can be viewed as a holomorphic 1-form with values in the adjoint bundle $Ad(E_G)= E_G\times_{Ad}\mathfrak{g}$.  If $P$ is any homogeneous polynomial on $\mathfrak{g}$ invariant under the conjugation action of $G$, then $P(\Phi)$ is a well defined section of $\Ksigma^r$ where $r$ is the degree of the polynomial. A basis for the ring of invariant polynomials thus defines a map, called the Hitchin fibration

\begin{equation}
\mathfrak{H}:\Mhiggs G)\rightarrow\bigoplus H^0(K^{r_i})
\end{equation}

\noi where the $r_i$ are the degrees of the invariant polynomials in the basis. Hitchin showed that the fibers of  this fibration are generically tori and that it gives $\Mhiggs G)$ the structure of a algebraic completely integrable system. This fibration has many ramifications; we focus on just one here, namely the existence of sections.

First constructed in \cite{liegroupsteichmuller},  Hitchin demonstrated the existence of  sections 

\begin{equation}\label{Hsection}
\Psi: \bigoplus H^0(K^{r_i})\rightarrow \Mhiggs G)
\end{equation}

\noi with the property that their image lies in $\Mhiggs G_{split})\subset \Mhiggs G)$, where $G_{split}$ denotes the split real form of the complex group $G$.  The construction is very explicit. The points in the image of $\Psi$ can be represented by Higgs bundles $(E,\Phi)$ in which the bundle $E$ is the same for all points of the section, and the Higgs field is described with reference to that bundle. For example, in the case of $G=\SL(2n,\C)$, the bundle in the image of the section can be chosen to be 

\begin{equation}\label{Hsectionbundle}
E=\Ksigma^{-n+\frac{1}{2}}\oplus \Ksigma^{-n+\frac{3}{2}}\oplus\cdots \Ksigma^{n-\frac{1}{2}}
\end{equation}

\noi and the Higgs fields are given by

\begin{equation}\label{companion}
\Phi=\begin{bmatrix}0&1&0&0&\cdots&0\\
\alpha_2&0&1&0&\cdots&0\\
\alpha_3&\alpha_2&0&1&\cdots&0\\
\vdots&\vdots&\vdots&\vdots&\cdots&\vdots\\
\alpha_{n-1}&\cdots&\alpha_3&\alpha_2&0&1\\
\alpha_n&\alpha_{n-1}&\cdots&\alpha_3&\alpha_2&0 \end{bmatrix}
\end{equation}

\noi Here, for appropriate $i\in\{2,3,\dots,n\}$ and $j\in\{\frac{1-2n}{2},\frac{3-2n}{2},\dots,\frac{2n-1}{2}\}$,
\begin{itemize}
\item $1$ is the constant unit length section in $H^0(\Ksigma^{-(j+1)}\otimes \Ksigma^{j}\otimes \Ksigma)=H^0(\mathcal{O}_{\Sigma})$,
\item  $\alpha_i\in H^0(\Ksigma^{-j}\otimes \Ksigma^{j+i-1}\otimes \Ksigma)=H^0(K^{i}_{\Sigma})$
\end{itemize}

\noi The Hitchin section in this case is defined by 

\begin{equation}
\Psi(\alpha_1,\dots,\alpha_n)=[(E,\Phi)]
\end{equation}

\noi where $(E,\Phi)$ are given by \eqref{Hsectionbundle} and \eqref{companion} respectively.  Hitchin proved in \cite{liegroupsteichmuller} that the resulting Higgs bundles are stable, i.e.\ that they define points in $\Mhiggs \SL(2n,\C))$. Moreover he showed that the image lies in $\Mhiggs \SL(2n,\R))$ and that the map defines a homeomorphism onto a connected component in $\Mhiggs \SL(2n,\R))$.

The construction generalizes so that for all complex semisimple $G$, the image of a Hitchin section is open and closed in $\Mhiggs G_{split})$ and thus defines a connected component.  The construction of these components presents them as complex vector spaces of complex dimension $(g-1)\dim(G_{split})$ parameterized by global holomorphic differentials. The precise descriptions are given in Table \ref{Hitchincomponents}. We will describe the construction in a more general context later (see the Cayley map in Section \ref{magicalSL2}) but it is instructive to point out some key elements in the case of $\SL(2n,\R)$. Foremost among those is the role of a principal $\mathfrak{sl}(2,\C)$ subalgebra in $\mathfrak{sl}(2n,\C)$, as described in the next section.

\subsection{The underlying $\mathfrak{sl}$-triple}

The structure group for the bundle in \eqref{Hsectionbundle} is $T=\C^*$. Let $\mathcal{E}_T$ be the frame bundle for $K^{\frac{1}{2}}_{\Sigma}$. Define the group homomorphisms

\begin{equation}\label{Tembed}
T\rightarrow\SL(2,\C)\rightarrow\SL(2n,\C)
\end{equation}

\noi where the first arrow is the map $t\mapsto\begin{bmatrix}t&0\\ 0&t^{-1}\end{bmatrix}$ and the second arrow is defined by the irreducible representation on $\C^{2n}$, i.e.\ defines a principally embedded $\mathfrak{sl}(2,\C)$ subalgebra in $\mathfrak{sl}(2n,\C)$.  Using \eqref{Tembed} to extend the structure group we can then construct the principal $\SL(2n,\C)$-bundle $\mathcal{E}_T[\SL(2n,\C)]$ thus realizing \eqref{Hsectionbundle} as the associated $\C^{2n}$-bundle. Moreover, viewed from this perspective the bundle of trace-free endomorphsims, i.e.\ $End_0(E)$, is the associated bundle

\begin{equation}
End_0(E)= \mathcal{E}_T[\SL(2n,\C)][\mathfrak{sl}(2n,\C)]
\end{equation}

Under the principal embedding, the image of the nilpotent element $f=\begin{bmatrix}0&1\\0&0\end{bmatrix}$ in $\mathfrak{sl}(2n,\C)$ is 

\begin{equation}\tilde{f}=\begin{bmatrix}0&1&0&0&\cdots&0\\
0&0&1&0&\cdots&0\\
0&0&0&1&\cdots&0\\
\vdots&\vdots&\vdots&\vdots&\cdots&\vdots\\
0&0&0&0&\cdots&1\\
0&0&0&0&\cdots&0\\ \end{bmatrix}
\end{equation}

Via \eqref{Tembed} the action of $T$ on $f$ and hence $\tilde{f}$ is given by $t\cdot \tilde{f}=t^{-2}\tilde{f}$. It follows that $\tilde{f}$ defines a global section of $End_0(E)\otimes\mathcal{E}_T^2=End_0(E)\otimes \Ksigma$.  

The adjoint action of the principal $\mathfrak{sl}(2,\C)$-subalgebra decomposes $\mathfrak{sl}(2n,\C)$ into irreducible $\mathfrak{sl}(2,\C)$-representations as

\begin{equation}
\mathfrak{sl}(2n,\C)=W_1\oplus W_2\oplus\dots\oplus W_{n-1}
\end{equation}

\noi where $\dim(W_m)=2m+1$ and $W_1$ is the principal $\mathfrak{sl}(2,\C)$. An element in $\mathfrak{sl}(2n,\C)$ of the form in \eqref{companion} (where now the $\alpha_i$ are just numbers) can then be expressed as

\begin{equation}\label{companionPhi}
\Phi=\tilde{f}+\alpha_1e_1+\dots\alpha_l e_l
\end{equation}

\noi where $e_i$ is a heighest weight vector for $W_i$.  The subspaces spanned by the $e_i$ (denoted by $<e_i>$) are preserved under conjugation by $T$, yielding subbundles 

\begin{equation}
\mathcal{E}_T[\SL(2n,\C)][<e_j>]\subset \mathcal{E}_T[\SL(2n,\C)][\mathfrak{sl}(2n,\C)]
\end{equation}

\noi Furthermore, if $\alpha_j$ is a section of $K^j_{\Sigma}$ then $\alpha_je_j$ is well defined as a section of  

\begin{equation}
\mathcal{E}_T[\SL(2n,\C)][<e_j>]\otimes \Ksigma\subset \mathcal{E}_T[\SL(2n,\C)][\mathfrak{sl}(2n,\C)]\otimes \Ksigma\ ,
\end{equation}

\noi thus recasting \eqref{companionPhi} as an expression of the form \eqref{companion}.

\begin{table}
\centering
  \resizebox{\textwidth}{!}{    
\begin{tabular}{|c|c|}
\hline
$G_{split}$ & Hitchin component\\
\hline
&\\ [-1em]
$\SL(n,\R)$ &$\bigoplus_{j=2}^n H^0(\Ksigma^j)$\\
&\\ [-1em]
\hline
&\\ [-1em]
$\SO(n,n+1)$ &$\bigoplus_{j=1}^n H^0(\Ksigma^{2j})$\\
&\\ [-1em]
\hline
&\\ [-1em]
$\mathrm{Sp}(2n,\R)$ &$\bigoplus_{j=1}^n H^0(\Ksigma^{2j})$\\
&\\ [-1em]
\hline
&\\ [-1em]
$\SO(n,n)$ &$\bigoplus_{j=1}^{n-1} H^0(\Ksigma^{2j})\oplus H^0(\Ksigma^n)$\\
&\\ [-1em]
\hline
&\\ [-1em]
$G_{2,split}$ &$H^0(\Ksigma^2)\oplus H^0(\Ksigma^6)$\\
&\\ [-1em]
\hline
&\\ [-1em]
$F_{4,split}$ &$H^0(\Ksigma^2)\oplus H^0(\Ksigma^6)\oplus H^0( \Ksigma^{8})\oplus H^0( \Ksigma^{12})$\\
&\\ [-1em]
\hline
&\\ [-1em]
$E_{6,split}$ &$H^0(\Ksigma^2)\oplus H^0(\Ksigma^5)\oplus H^0( \Ksigma^{6})\oplus H^0( \Ksigma^{8})\oplus H^0( \Ksigma^{9})\oplus H^0( \Ksigma^{12})$\\
&\\ [-1em]
\hline
&\\ [-1em]
$E_{7,split}$ &$H^0(\Ksigma^2)\oplus H^0(\Ksigma^6)\oplus H^0( \Ksigma^{8})\oplus H^0( \Ksigma^{10})\oplus H^0( \Ksigma^{12})\oplus H^0( \Ksigma^{14})\oplus H^0( \Ksigma^{18})$\\
&\\ [-1em]
\hline
&\\ [-1em]
$E_{8,split}$ &$H^0(\Ksigma^2)\oplus H^0(\Ksigma^8)\oplus H^0( \Ksigma^{12})\oplus H^0( \Ksigma^{14})\oplus H^0( \Ksigma^{18})\oplus H^0( \Ksigma^{20})\oplus H^0( \Ksigma^{24})\oplus H^0( \Ksigma^{30})$\\
\hline
\end{tabular}
}
\medskip
\caption{Hitchin/Teichm\"uller  components in $\Mhiggs G_{split})$}\label{Hitchincomponents}
\end{table}

\begin{remark} Hitchin points out (Section \S4 in \cite{liegroupsteichmuller}) that Theorem 7 of \cite{K} proves that there exists a basis $p_1,\dots p_l$ of invariant polynomials on $\mathfrak{g}$ such that

\begin{equation}
p_i(f+\alpha_1e_1+\dots\alpha_l e_l)=\alpha_i
\end{equation}

\noi where $e_i$ are heighest weight vectors.  Moreover, the degree of $p_i$ is $2m_i+1$, where $2m_i+1$ is the dimension of $V_i$ or equivalently $m_i$ is the eigenvalue of $ad_h$ on $e_i$. These fundamental invariants of the Lie algebra are called the exponents.
\end{remark}

\subsection{Local minima of the Hitchin function}

The components of the moduli space $\Mhiggs G_{split})$  constructed via sections of the Hitchin fibration are called Hitchin or Teichm\"uller  components. The Hitchin function can be evaluated on these components.  It follows from the description of the Higgs fields as in \eqref{companionPhi} that the local minimum occurs at the Higgs bundle where the $\alpha_i$ are all zero, i.e.\ at the image of zero in the Hitchin section. The Higgs field at this minimum is $\Phi=\tilde{f}$. In particular it is not zero; indeed the Higgs field vanished nowhere on Hitchin components. The components constructed in this way are thus clearly distinct from components with trivial local minima.

The Hitchin components  are examples of higher Teichm\"uller  spaces  (see for example \cite{W2018} for more about this).  Another source of higher Teichm\"uller  spaces comes from real forms of Hermitian type, where distinguished components arise in a seemingly different way, as described in the next section.   In Section \ref{magicalSL2} we will see that the underlying mechanism is actually the same as for the Hitchin components, with special $\mathfrak{sl}(2,\C)$-subalgebras providing the key.

\section{Maximal components and the Cayley correspondence}\label{maxcomponents}

Suppose that $G$ is a non-compact real form with maximal compact subgroup $H$ such that the symmetric space $G/H$ admits a $G$-invariant complex structure. Such real forms are said to be of Hermitian type. The following\footnotemark\footnotetext{Using Helgason's notation as in \cite{Helgason}} is a complete list of the corresponding real Lie algebras if $G$ is simple and connected:

$$\mathfrak{su}(p,q),\ \mathfrak{sp}(2n,\R),\ \mathfrak{so}(2,q),\ \mathfrak{so}^*(2n),\ \mathfrak{e}_6^{-14},\ \mathrm{and}\ \mathfrak{e}_7^{-25}\ . $$

These fall into two classes, referred to as {\it tube type} or {\it non-tube type}. The non-tube  cases are $\mathfrak{su}(p,q)$ with $p\ne q$, $\mathfrak{so}^*(2n)$ with $n$ odd, and $\mathfrak{e}_6^{-14}$; the rest are called tube type because their associated symmetric spaces $G/H$ are biholomorphic to tubes, i.e.\ domains $V+i\Omega$ where $V$ is a real vector space and $\Omega\subset V$ is a proper convex open cone. 

Higgs bundles for real forms of Hermitian type have been studied on a case-by-case basis (see the references in the Tables) and more systematically in \cite{BGRmaximalToledo}. All the results are now subsumed under the umbrella of the Cayley components determined by magical $\mathfrak{sl}_2$-triples, as described in Section \ref{magicalSL2}.

For all real form $G$ of Hermitian type a discrete invariant called the Toledo number can be attached to a $G$-Higgs bundle or equivalently (via the NAH correspondence) to a surface group representation into $G$. The invariant is subject to a bound known as a Milnor-Wood bound. This allows one to distinguish {\it maximal} Higgs bundles or representations as those for which the Toledo invariant attains its maximal value.  Since the invariant is constant on connected components, this leads to the following definition

\begin{definition} Let $G$ be a real form of Hermitian type. A component of $\Mhiggs G)$ (respectively $\mathrm{Rep}(\pi_1(S),G)$) is called a {\it maximal component} if the Toledo invariant attains its maximal value on the Higgs bundles (respectively surface group representations) represented in that component. 
\end{definition}

If $G$ is of tube type then the maximal components in $\mathrm{Rep}(\pi_1(S), G)$ have several special features, by virtue of which they are considered examples of higher Teichm\"uller  spaces (see \cite{W2018} for an excellent overview).

In the Higgs bundle moduli space, the special nature of the maximal components manifests as extra symmetries on the bundles and Higgs fields.  As a result of these constraints the maximal components of $\Mhiggs G)$ acquire new descriptions as moduli spaces of Higgs-like pairs but for a new group, $G'$, called the Cayley partner to $G$.  We illustrate this phenomenon in the case of $G=\mathrm{Sp}(2n,\R)$.

\subsection{$\mathrm{Sp}(2n,\R)$-Higgs bundles}\hfill

The maximal compact subgroup of $\mathrm{Sp}(2n,\R)$ is $H=U(n)$ and the complexified Cartan decomposition can be given as

\begin{align}
\mathfrak{sp}(2n,\C)=&\mathfrak{h}^{\C}\oplus \mathfrak{m}^{\C},\ \mathrm{where}\label{h+m}\\
\mathfrak{h}^{\C}=&\Big\{\begin{bmatrix}A&0\\0&-A^T\end{bmatrix}\ \Big |\  A\in\mathfrak{gl}(n)\Big\}\ ,\\
\mathfrak{m}^{\C}=&\Big\{\begin{bmatrix}0&C\\D&0\end{bmatrix}\ \Big |\  C,D\in\mathfrak{gl}(n)\ , C=C^T, D=D^T\Big\}\label{mc}
\end{align}

\noi An $\mathrm{Sp}(2n,\R)$-Higgs bundle can thus be described as a rank $2n$ Higgs vector bundle of the form $(V\oplus V^*, \Phi)$ where
\begin{itemize}
\item $V$ is a rank $n$ holomorphic vector bundles, and
\item $\Phi=\begin{bmatrix}0&\beta\\\gamma&0\end{bmatrix}$ with $\beta\in H^0(Sym^2(V)\otimes \Ksigma)$ and $\gamma\in H^0(Sym^2(V^*)\otimes \Ksigma)$.
\end{itemize}

\noi The Toledo invariant is defined to be

\begin{equation}
\tau(V\oplus V^*, \Phi)=\deg(V)
\end{equation}

\noi and satisfies the bound

\begin{equation}
|\tau(V\oplus W, \Phi)|\le n(g-1)
\end{equation}

\begin{theorem} If $(V\oplus V^*, \Phi)$ is a semistable $\mathrm{Sp}(2n,\R)$-Higgs bundle with maximal Toledo invariant $\tau= n(g-1)$ then $\gamma$ defines an isomorphism $\gamma:V\simeq V^*\otimes \Ksigma$ with $\gamma\in H^0(Sym^2(V^*)\otimes \Ksigma)$.\end{theorem}

\noi Thus if we fix a square root of $\Ksigma$, denoted by $K^{-1/2}_{\Sigma}$, and define $U=V\otimes K^{-\frac{1}{2}}_{\Sigma}$, then 

\begin{itemize}
\item $q_U=\gamma\otimes 1_{K^{-\frac{1}{2}}_{\Sigma}}:U\rightarrow U^*$ defines an orthogonal structure, and
\item $\theta=\beta\circ\gamma:U\rightarrow U\otimes \Ksigma^2$ defines a $\Ksigma^2$-twisted symmetric endomorphism.
\end{itemize}

\noi Notice that the structure group for the orthogonal bundle is $O(n,\C)$ and the $\Ksigma^2$-twisted Higgs field takes values in $Sym^2(\mathfrak{gl}(n,\C))$.  Moreover, the complexified Cartan decomposition for $\mathfrak{gl}(n,\R)$ is

\begin{equation}
\mathfrak{gl}(n,\C)=\mathfrak{o}(n,\C)\oplus Sym^2(\mathfrak{gl}(n,\C)
\end{equation}

\noi Thus the defining data for maximal $\mathrm{Sp}(2n,\R)$-Higgs bundles can be identified with the defining data for $\Ksigma^2$-twisted $G'$-Higgs bundle where $G'=GL(n,\R)$. 

\begin{theorem}The Cayley partner for $\mathrm{Sp}(2n,\R)$ is $GL(n,\R)$ and  the map

\begin{equation}\label{cayleymap}
(U,\theta)\mapsto \Big ((U\otimes K^{\frac{1}{2}}_{\Sigma})\oplus (U^*\otimes K^{-\frac{1}{2}}_{\Sigma}),\begin{bmatrix}0&\theta\circ q_U^{-1}\otimes 1\\ q_u\otimes 1&0  \end{bmatrix}\Big )
\end{equation}

defines an isomorphism (the Cayley map)

\begin{equation}
\Psi:\mathcal{M}_{K^2}(\Sigma,GL(n,\R))\rightarrow\mathcal{M}_{max}(\Sigma,\mathrm{Sp}(2n,\R))
\end{equation}

\end{theorem}

If $G$ is any other simple group of Hermitian tube type, the maximal component $\mathcal{M}_{max}(\Sigma,G)$ similarly
acquires a description as a moduli space of $K^2$-twisted $G'$-Higgs bundles where $G'$ is the Cayley partner group.  A list of the Cayley partners can be found in \cite{BGRmaximalToledo}, where it is also shown that the Cayley partner of $G$ can be characterized as the non-compact dual of $H$ in $H^{\C}$, where $H\subset G$ is a maximal compact subgroup and $H^{\C}$ is its complexification.  

The Cayley map gives a description of the maximal components which, while not as explicit as the parametrization of the Hitchin components, can nevertheless be useful in understanding global features of the component.

Like the special features of the Hitchin section, the appearance of these Cayley partners can be traced back to a distinguished $\mathfrak{sl}_2$-triple in the complexified Lie algebra. Indeed we will see a unifying construction in Section \ref{magicalSL2} but we briefly describe the specific case of $\mathrm{Sp}(2n,\R)$ to illustrate the way it works in the case of real forms of Hermitian tube type.

\subsection{The underlying $\mathfrak{sl}_2$-triple} \hfill

As above, if $(V\oplus V^*,\Phi)$ is a polystable $\mathrm{Sp}(2n,\R)$-Higgs bundle with maximal Toledo invariant then 

\begin{equation} 
V\oplus V^*\ =\ \big(U\otimes K^{\frac{1}{2}}_{\Sigma}\big)\oplus \big(U^*\otimes K^{-\frac{1}{2}}_{\Sigma}\big)
\simeq U\otimes \big(K^{\frac{1}{2}}_{\Sigma}\oplus K^{-\frac{1}{2}}_{\Sigma}\big)
\end{equation}

The frame bundle for  $V\oplus V^*$ is thus of the form $\mathcal{E}_{O(n,\C)}\star\mathcal{E}_{T}$ where 

\begin{itemize}
\item $\mathcal{E}_{O(n,\C)}$ is the principal $O(n,\C)$-bundle which defines $U$,
\item $\mathcal{E}_{T}$ is the principal $\C^*$-bundle which defines $K^{\frac{1}{2}}_{\Sigma}$, and
\item the $\star$ operation is a multiplication in $\mathrm{Sp}(2n,\C)$, using the embeddings
\begin{align}
A&\mapsto \begin{bmatrix}A&0\\0&A\end{bmatrix}\ \mathrm{for}\ A\in O(n,\C),\label{On}\\
\lambda&\mapsto \begin{bmatrix}\lambda I&0\\0&\lambda^{-1}I\end{bmatrix}  \mathrm{for}\ \lambda\in\C^*
\end{align}
 \end{itemize}

\noi Notice that

\begin{itemize}
\item the $\C^*$ in $\mathrm{Sp}(2n,\C)$ lies in a copy of $\SL(2,\C)$ embedded as $S=\SL(2,\C)\otimes I_n$, i.e.\ via

\begin{equation}
\begin{bmatrix} a&b\\c&d\end{bmatrix}\mapsto \begin{bmatrix} aI_n&bI_n\\cI_n&dI_n\end{bmatrix}
\end{equation}
\noi Moreover, within $\Sp(2n,\C)$  this $\SL(2,\C)$-subgroup  commutes with the $O(n,\C)$-subgroup defined by \eqref{On}, thereby allowing $\mathcal{E}_{O(n,\C)}\star\mathcal{E}_{T}$ to be well defined as an $\Sp(2n,\C)$-bundle.

\item the summand $\mathfrak{m}^{\C}$ in \eqref{mc} decomposes as $\mathfrak{m}^{\C}_+\oplus\mathfrak{m^{\C}}_-$ where $D=0$ in $\mathfrak{m}^{\C}_+$ and $C=0$ in $\mathfrak{m}^{\C}_-$.  Moreover if 

\begin{equation}T=\Big\{\begin{bmatrix}t I&0\\0&t^{-1}I\end{bmatrix}\ \Big |\ t\in\C\Big\}
\end{equation}

\noi then the adjoint action of $T$ has weight $0$ on $\mathfrak{h}^{\C}$ and weights $\pm 2$ on $\mathfrak{m}_{\pm}^{\C}$. The decomposition in \eqref{h+m} can thus be identified as

\begin{align}
\mathfrak{sp}(2n,\C)&=\mathfrak{m}^{\C}_-\oplus\mathfrak{h}^{\C}\oplus \mathfrak{m}^{\C}_+\\
&=\mathfrak{sp}(2n,\C)_{-2}\oplus \mathfrak{sp}(2n,\C)_0\oplus\mathfrak{sp}(2n,\C)_2\label{C*forsp2n}
\end{align}

\noi where the subscripts denote the weights of the $T$ action

\item The Lie algebra of $S$ is an $\mathfrak{sl}(2,\C)$-subalgebra spanned by the triple

\begin{equation}
\{f,h,e\}=\Big\{\begin{bmatrix} 0&0\\ I_n&0\end{bmatrix}, \begin{bmatrix} I_n&0\\ 0&-I_n\end{bmatrix}, \begin{bmatrix} 0&I_n\\ 0&0\end{bmatrix}\Big\}
\end{equation}

\noi Under the adjoint action the decomposition of $\mathrm{Sp}(2n,\C)$ into irreducible $\mathfrak{sl}(2,\C)$-representations takes the form

\begin{equation}\label{sl2forsp2n}
\mathfrak{sp}(2n,\C)=W_0\oplus W_1
\end{equation}

\noi where $W_m$ is a sum of irreducible representations of dimension $2m+1$. 
\item The decompositions in \eqref{sl2forsp2n} and \eqref{C*forsp2n} are related by the fact that $\mathfrak{sp}(2n,\C)_2$ is the highest weight space for $W_1$.
\end{itemize}

\noi The Higgs field can thus be decomposed as 

\begin{equation}
\Phi=\tilde{f}+\tilde{\theta}
\end{equation}

\noi  where $\tilde{f}$ is defined by the nilpotent $f$ in $\mathfrak{sl}(2,\C)$ and $\tilde\theta=\begin{bmatrix}0&\theta\circ q_U^{-1}\otimes 1\\0&0  \end{bmatrix}$ is determined by a highest weight element.
 
\section{Magical triples and generalized Cayley maps}\label{magicalSL2}

The underlying mechanism responsible for the Hitchin components when $G$ is a split real form, or the maximal components when $G$ is Hermitian of tube type, can be understood in a uniform way by focussing on the role of the distinguished $\mathfrak{sl}_2$-triples in the Lie algebras. From this perspective, the Hitchin components and the maximal components represent opposite extremes in a range which includes special components in the moduli spaces for $G$-Higgs bundles for $G$ on a list which also includes $\SO(p,q)$ or certain exceptional real groups.  The list of real forms which occur in this range is the same as the list which admit a positivity structure known as $\Theta$-positivity. Introduced in \cite{GW18}, the $\Theta$-positivity property plays a crucial role in higher Teichm\"uller  theory, where it is responsible for many of the similarities between Teichm\"uller  space and special components of $\mathrm{Rep}(\pi_1(S),G)$ known collectively as higher Teichm\"uller  spaces (see \cite{W2018,GW22})

\subsection{Magical triples}\hfill

Let $\mathfrak{g}^{\C}$ be a complex semisimple Lie algebra.  By the Jacobson-Morozov theorem any non-zero nilpotent element $e\in\mathfrak{g}^{\C}$ can be completed to a triple (called an $\mathfrak{sl}_2$-triple) of nonzero elements $\{f, h, e\}$ satisfying

\begin{equation}
[h,e]=2e\ ,\ [h,f]=-2f\ ,\ [e,f]=h\ ,
\end{equation}

\noi i.e.\ which generate an $\mathfrak{sl}(2,\C)$-subalgebra with $h$ semisimple.  Any $\mathfrak{sl}_2$-triple in $\mathfrak{g}^{\C}$ determines two decompositions of  $\mathfrak{g}^{\C}$, namely

\begin{enumerate}
\item into weight-spaces for the ad-action of $h$, denoted by 
\begin{equation}\label{weights}
\mathfrak{g}^{\C}= \mathfrak{g}_0\oplus\bigoplus_{j=1}^M\mathfrak{g}_j
\end{equation}
\noi where the subscript denotes the weight, and
\item into irreducible $\mathfrak{sl}(2,\C)$-representations, denoted
\begin{equation}\label{sl2reps}
\mathfrak{g}^{\C}= W_0\oplus\bigoplus_{j=1}^MW_{m_j}
\end{equation}
\noi where $W_{m_j}$ is isomorphic to a direct sum of $n_j$ copies (with $n_j> 0$) of the unique $(m_j + 1)$- dimensional representation.  The integer pairs $\{(m_j,n_j)\ | j=1,2,\cdots,M\}$ from \eqref{weights} and \eqref{sl2reps} are called the {\it $\mathfrak{sl}_2$-data for the triple $\{f,h,e\}$}.
\end{enumerate}

\noi Note: The subspace of 1-dimensional representations, i.e.\ the summand $W_0$ in \eqref{sl2reps}, is the centralizer of the the $\mathfrak{sl}(2,\C)$-subalgebra generated by the $\mathfrak{sl}_2$-triple. We denote this by $\mathfrak{c}$.

A nilpotent $e$ is called {\it even} if all weights in \eqref{weights} are even. Since the highest $ad_h$-weight in $W_{k}$ is $k$, for $\mathfrak{sl}(2,\C)$-subalgebras determined by even nilpotents we get decompositions

\begin{equation}
W_{2m_j}=V_{2m_j}\oplus ad_f(V_{2m_j})\oplus ad^2_f(V_{2m_j})\oplus\cdots\oplus ad_j^{2m_j}(V_{2m_j})
\end{equation}

\noi where $V_{2m_j}=W_{2m_j}\cap\mathfrak{g}_{2m_j}$ is the heighest weight space and $Z_{2m_j}=ad_f^{m_j}(V_{2m_j})\subset \mathfrak{g}_0$.

\begin{definition}Define a vector space involution $\sigma:\mathfrak{g}^{\C}\rightarrow\mathfrak{g}^{\C}$ by the linear extension of 

\begin{equation}\label{sigma}
\sigma (x)=\begin{cases} x\ \mathrm{if}\ x\in  V_0\\
(-1)^{k+1}  \mathrm{if}\ x\in  ad_f^k(V_j)\ \mathrm{for\ some}\ 0\le k\le j\ \mathrm{and}\  j>0
\end{cases}
\end{equation}
\end{definition}

\noi In general $\sigma$ does not preserve the Lie algebra structure on $\mathfrak{g}^{\C}$.  In \cite{BCGGO21} we introduced the notion of a magical {\it $\mathfrak{sl}_2$-triple} defined as follows.

\begin{definition} The $\mathfrak{sl}_2$-triple $\{f,h,e\}$ is called {\it magical} if the involution $\sigma$  preserves the Lie algebra structure on $\mathfrak{g}^{\C}$.
\end{definition}

We have already seen two examples of magical $\mathfrak{sl}_2$-triples, namely
\begin{itemize}
\item principal $\mathfrak{sl}_2$-triples, and
\item the triples associated to real forms of Hermitian tube type.
\end{itemize}

If the involution $\sigma$ is magical, and thus a Lie algebra involution, then it has an associated anti-holomorphic involution\footnotemark\footnotetext{The antiholomorphic involution, say $\tau$ is unique up to conjugation and satisfies the condition that $\sigma\tau=\tau\sigma$ defines a compact real form.}
 which defines a real form for $\mathfrak{g}^{\C}$. 
 
\begin{definition} If $\{f,h,e\}$ is a magical $\mathfrak{sl}_2$-triple and $\tau_e$ is the antiholomorphic involution associated to the involution defined by \eqref{sigma}, the real form of $\mathfrak{g}^{\C}$ defined by $\tau_e$ is called the {\it canonical real form associated to the magical $\mathfrak{sl}_2$-triple, or simply the magical real form}.

Let $G^{\C}$ be a connected complex Lie group with Lie algebra $\mathfrak{g}^{\C}$ and such that the involution $\sigma$ integrates to an involution on $G^{\C}$.  Denote also by $\tau_e:G^{\C}\rightarrow G^{\C}$ the involution which integrates $\tau_e$. The fixed-point group defined by $\tau_e$ is called the {\it canonical real form of $G$ defined by the magical $\mathfrak{sl}_2$-triple, or simply the magical real form of $G^{\C}$}. 
\end{definition} 
 
 \begin{theorem} [\cite{BCGGO21}] \label{classification}  A real form of a complex simple Lie algebra is the canonical real form associated to the magical $\mathfrak{sl}_2$-triple if and only if it is one of the following types:
\begin{enumerate}
\item a split real form,
\item a real form of Hermitian tube type,
\item $\mathfrak{so}(p,q)$ with $2<p< q$,
\item the quaternionic real form of $\mathfrak{f}_4, \mathfrak{e}_6, \mathfrak{e}_7$, or $\mathfrak{e}_8$ (a.k.a. $\mathfrak{f}_4^4,\mathfrak{e}_6^2,\mathfrak{e}_7^{-5}$, or $\mathfrak{e}_8^{-24}$)
\end{enumerate} 
\end{theorem}

\begin{remark} The four types in Theorem \eqref{classification} are not mutually exclusive. The overlaps occur because the following real forms are split and hence also included in type (1):
\begin{itemize}
\item $\mathfrak{sp}(2n,\R)$ (in type (2))
\item $\mathfrak{so}(n,n+1)$  (in type (3))
\item $\mathfrak{f}_4^4$ (in type (4))
\end{itemize}
Each of these real forms admits two distinct magical $\mathfrak{sl}_2$-triples, with the characteristics of each $\mathfrak{sl}_2$-triple determined by one of the two types to which the real form belongs.

\end{remark}

\subsection{Structures determined by magical $\mathfrak{sl}_2$-triples}\hfill

In the next three subsections we assume that $\{f,h,e\}$ is a magical subtriple in $\mathfrak{g}^{\C}$, and that $G^{\C}$ is as above, i.e.\ a connected complex Lie group with Lie algebra $\mathfrak{g}^{\C}$ and such that the involution $\sigma$ integrates to an involution on $G^{\C}$.  We note that in this case $e$ is an even nilpotent so the subspaces $Z_{2j}=W_{2j}\cap\mathfrak{g}_0$ are non-empty.  We denote the canonical real form associated to the magical $\mathfrak{sl}_2$-triple by $G$ and its maximal compact subgroup by $H$.  Then in addition to the decompositions \eqref{weights} and \eqref{sl2reps} we have the Cartan decomposition

\begin{equation}\label{magical cartan}
\mathfrak{g}^{\C}=\mathfrak{h}^{\C}\oplus\mathfrak{m}^{\C}
\end{equation}

\noi where $\mathfrak{h}^{\C}$ is the complexification of $\mathfrak{h}=Lie(H)$.   We get the following structures (see \cite{BCGGO21} for details):
\begin{enumerate}
\item {\bf A real form of $\mathfrak{g}_0$:}  The summand $\mathfrak{g}_0$ in \eqref{weights}, i.e.\ in the zero weight-space for $ad_h$, is a subalgebra of $\mathfrak{g}^{\C}$. If $h$ is part of a magical $\mathfrak{sl}_2$-triple then $\mathfrak{g}_0$ admits a real form - called the Cayley real form - defined by the involution

\begin{align}
&\theta: \mathfrak{g}_0\rightarrow \mathfrak{g}_0\\
&\theta= \begin{cases}
 +1\ \mathrm{on}\ W_0\\ -1\ \mathrm{on}\ W_j\cap\mathfrak{g}_0\ \mathrm{for}\ j > 0
 \end{cases}
\end{align}

\noi Denote this real form by $\mathfrak{g}_C$.  

\item {\bf A decomposition of $\mathfrak{g}_C$:}  If the magical $\mathfrak{sl}_2$-triple defined by the nilpotent $e$ is not principal then there is a subalgebra $\mathfrak{g}(e)\subset \mathfrak{g}$ such that the $\mathfrak{sl}_2$-triple is principal in $\mathfrak{g}(e)$. In fact $\mathfrak{g}(e)$ is the centralizer of $\mathfrak{c}=W_0$ as defined in \eqref{sl2reps}. The real form $\mathfrak{g}_C$ decomposes as

\begin{equation}
\mathfrak{g}_C = \mathfrak{g}_{0,ss}\oplus\R^{rank(\mathfrak{g}(e))}
\end{equation}

\noi where $\mathfrak{g}^{\R}_{0,ss}$ is either zero or a simple real Lie algebra.  
 
\item {\bf The Cayley group:}  The Lie algebra $\mathfrak{g}_C$ is does not uniquely determine a Cayley Lie group because it is insensitive to the center of the group. The requisite extra information is obtained as follows.  Let $S\subset G^{\C}$ be the subgroup whose Lie algebra is the magical $\mathfrak{sl}_2$-triple. In the adjoint form of $\mathfrak{g}^{\C}$ this is isomorphic to $\PSL(2,\C)\simeq\SO(3,\C)$; otherwise for the groups we consider, i.e.\ complex simple Lie groups, it isomorphic to $\SL(2,\C)$.  Let $C$ be the centralizer of $S$ in $G^{\C}$.   
\end{enumerate}
\begin{definition} The Cayley group of the magical $\mathfrak{sl}_2$-triple in $G^{\C}$ is  the group

\begin{equation}\label{cayleygroup intro}
\mathrm{G}_C=\mathrm{G}_{0,ss}\times (\R^+)^{rank(\mathfrak{g}(e))},
\end{equation}

\noi where $\mathrm{G}_{0,ss}\subset G^{\C}$ is a real Lie group with Lie algebra $\mathfrak{g}_{0,ss}$ and maximal compact subgroup $C\cap G$
\end{definition}

\subsection{Special $G$-Higgs bundles}\hfill

Recall that a $G$-Higgs bundle is a pair $(E_{H^{\C}},\Phi)$ (see Section \ref{GHiggs}).  If $G$ is a magical real form then the structures described above can be used to construct $G$-Higgs bundles with special structure. These turn out to comprise a union of connected components in $\Mhiggs G)$.

By definition, the subgroups $S$ and $C$ commute in $G^{\C}$. It follows that principal $S$- and $C$-bundles, say $\mathcal{E}_S$ and $\mathcal{E}_C$, can be combined to produce a principal $G^{\C}$-bundle whose structure group factors through the image of $S\times C$ in $G^{\C}$. We denote this construction by $\mathcal{E}_S\star\mathcal{E}_C$.  In particular, if $T\subset S$ is the $\C^*$-subgroup whose Lie algebra is spanned by $h$, and $\mathcal{E}_T$ is a principal $T$-bundle, then we can construct $\mathcal{E}_T\star\mathcal{E}_C$. We will take $\mathcal{E}_T$ to be the frame bundle for $\Ksigma$ if $S\simeq \PSL(2,\C)$, or a square root $K^{\frac{1}{2}}_{\Sigma}$ if $S\simeq \SL(2,\C)$.

A important subtlety is the difference between $C$ and $C\cap H^{\C}$, where $H^{\C}$ is the complexification of a maximal compact subgroup of $G$. From \eqref{cayleygroup intro} it is clear that $C\cap H^{\C}$ is the complexification of a maximal compact subgroup of the Cayley group $G_C$. Moreover, $C$ and $C\cap H^{\C}$ have the same Lie algebra, namely $W_0$. They may however have different centers. For example, for $G_{\R}=\SL(n,\R)$ we get $C=\Z_n$ but $H=\SO(n,\C)$ and $C\cap \SO(n,\C)$ is either $\Z_2$ ($n$ even) or trivial ($n$ odd) [See Remark 4.12 in \cite{BCGGO21}]

\begin{theorem}Suppose that the Lie algebra of $S$ is magical. Then $\mathcal{E}_T\star\mathcal{E}_C$ defines a principal $H^{\C}$-bundle if and only if the structure group of $\mathcal{E}_C$ reduces to $C\cap H^{\C}$. 
\end{theorem}

If  $\mathcal{E}_T\star\mathcal{E}_C$ defines a principal $H^{\C}$-bundle 

\begin{equation}
\mathcal{E}_{H^{\C}}=\mathcal{E}_T\star\mathcal{E}_C[H^{\C}]
\end{equation}

\noi then this can be used to construct $G$-Higgs bundles.  In particular, recall that the maximal compact subgroup of $G_{0,ss}$ can be taken to be $H_{0,ss}=C\cap G$ (see \eqref{cayleygroup intro}). Hence if $\mathcal{E}_{C}$ is a principal $H_{0,ss}^{\C}$-bundle (where $H_{0,ss}^{\C}$ denotes the complexification), then its structure group does reduce to $H^{\C}$. We assume this to be the case.

The Higgs fields must take values in summand $\mathfrak{m}^{\C}$ in \eqref{magical cartan}, i.e.\ in the $-1$-eigenspace for the involution defined by \eqref{sigma}. We get contributions from two sources:

\begin{enumerate}
\item Under the ad-action $S$ and $C$ preserve the highest weights subspaces $V_{2m_j}\subset W_{2m_j}$ in \eqref{sl2reps}, and by construction $V_{2m_j}\subset \mathfrak{m}^{\C}$.  We thus get subbundles 

\begin{equation}
\mathcal{E}_{H^{\C}}[V_{2j}]\subset \mathcal{E}_{H^{\C}}[\mathfrak{m}^{\C}]
\end{equation}

\item Since $C$ acts trivially on $f$ and the subspace spanned by $f$ is preserved by $T$, we get  a subbundle  

\begin{equation}
\mathcal{E}_{H^{\C}}[<f>]=\mathcal{E}_T[<f>]\subset \mathcal{E}_{H^{\C}}[\mathfrak{m}^{\C}]
\end{equation}

\end{enumerate}

\noi Recall that we take $\mathcal{E}_T$ to be the principal bundle for $\Ksigma$ or $K^{\frac{1}{2}}_{\Sigma}$, depending on whether $S$ is  $S=\SL(2,\C)$ or $\PSL(2,\C)$. Since $\C^*=T$ acts with weight $-2$ on $f$ if $S=\SL(2,\C)$ and weight $-1$ if $S=\PSL(2,\C)$, in both cases we get 

\begin{equation}\mathcal{E}_T[<f>]\otimes \Ksigma\simeq\mathcal{O}
\end{equation}

\noi and hence  $f\in \mathfrak{m}^{\C}$ defines a constant section, which we denote also by $f$.  In both cases we will refer to the pair $(\mathcal{E}_T,f)$ as the {\it uniformizing Higgs pair}.

Given any sections $\phi_{m_j}\in H^0(\mathcal{E}_{H^{\C}}[V_{2m_j}]\otimes \Ksigma$ we can thus construct $G$-Higgs bundles of the special form

\begin{equation}\label{specialHiggs}
\big (\mathcal{E}_T\star\mathcal{E}_C[H^{\C}], f+\phi_{m_1}+\cdots +\phi_{m_M}\big )
\end{equation}

\subsection{Cayley map and exotic components}\hfill

The section $f$ allows us to globalize the $ad_f$-action on $\mathfrak{g}^{\C}$ to define isomorphisms

\begin{equation}
ad_f^{m_j}: \mathcal{E}_{H^{\C}}[V_{2m_j}]\otimes \Ksigma\rightarrow \mathcal{E}_C[Z_{2m_j}]\otimes \Ksigma^{1+j}
\end{equation}

\noi where as above $Z_{2m_j}=W_{2m_j}\cap\mathfrak{g}_0$.

Two key facts are important:

\begin{enumerate}
\item If $Z_{2m_j}\simeq\C$ (i.e.\ $\dim(Z_{2m_j})=1$) then $\mathcal{E}_{H^{\C}}[Z_{2m_j}]\simeq\mathcal{O}$ and we can identify $ad_f^j(\phi_{m_j})$ with a section of $\Ksigma^{1+m_j}$, i.e.\ 

\begin{equation}ad_f^{m_j}(\phi_{m_j})=q_{m_j}\in H^0(\Ksigma^{1+m_j})\ .
\end{equation}
\item If the $\mathfrak{sl}_2$-triple is magical then $\dim(Z_{2m_j})>1$ for at most one $j>0$. Indeed, in the four cases listed in Theorem \ref{classification} this is trivially true in cases (1) and (2) ($\dim(Z_{2m_j})=1$ for all $j$ in case (1), and in case (2) there is exactly one $j>0$) and follows from calculations in the other cases (see theorem 4.3 in \cite{BCGGO21}).

\end{enumerate}

Let $m_c$ be the unique positive integer such that $\dim(Z_{2m_c})>1$. Recall that by definition, $Z_{2m_j}$ lies in the $-1$-eigenspace for the involution that defines the Cayley real form. If $\phi_{m_c}$ is a section of  $\mathcal{E}_C[V_{2{m_c}}]\otimes \Ksigma$ and $\psi_{m_c}=ad_f^{m_c}(\phi_{m_c})\in  \mathcal{E}_{H^{\C}}[Z_{2j}]\otimes \Ksigma^{2m_c+1}$ it then follows that $(\mathcal{E}_C,\psi_{m_c})$ defines a $\Ksigma^{2m_c+1}$-twisted $G_{0,ss}$-Higgs bundle. The remaining indices $m_j$, i.e.\ those for which $\dim(Z_{2m_c})=1$ are determined by the exponents for the subalgebra $\mathfrak{g}(e)$ (see Lemma 5.7 in \cite{BCGGO21}).

The $G$-Higgs bundles of the form in \eqref{specialHiggs} are thus determined by
\begin{itemize}
\item (Fixed data) $(\mathcal{E}_T,f)$,  i.e.\ the uniformizing Higgs pair, and
\item (Cayley parameters) $\big ((\mathcal{E}_C,\psi), q_1,\cdots, q_N)\big )$, i.e.\ a $\Ksigma^{2m_c+1}$-twisted $G_{0,ss}$-Higgs bundle together with a collection of $N$ holomorphic differentials where $N=rank(\mathfrak{g}(e)$).
\end{itemize}

The $G$-Higgs bundles constructed in this way are polystable if and only if the $\Ksigma^{2m_c+1}$-twisted $G_{0,ss}$-Higgs bundle are polystable. We thus get a map 

\begin{equation}\label{eq:Cayleymoduli}
\Psi_e:\mathcal{M}_{K^{m_c+1}}(G_{0,ss}^\R)\times \bigoplus_{j=1}^{\mathrm{rank}(\mathfrak{g}(e))}H^0(K^{l_j+1}_{\Sigma})\longrightarrow\mathcal{M}(G).
\end{equation}

\noi where the exponents $l_j$ are the exponents of the Lie algebra $\mathfrak{g}(e)$.

\begin{theorem}  [Theorem 7.1 and Corollary 7.6 in \cite{BCGGO21}] The image of the map \eqref{eq:Cayleymoduli}, which we call the {\it Cayley map}, is open and closed, and hence defines a union of connected component of $\mathcal{M}(G)$ homeomorphic to 

\begin{equation}\label{Cayleyproduct}
\mathcal{M}_{K^{m_c+1}}(G_{0,ss}^\R)\times \bigoplus_{j=1}^{\mathrm{rank}(\mathfrak{g}(e))}H^0(K^{l_j+1}_{\Sigma})\ .
\end{equation}

\noi Each such component is locally irreducible and irreducible.
\end{theorem}

\begin{table}\label{MagicalLiedata}
\begin{tabular}{|c||c|c|c|c|}
\hline
&&&&\\
Real form& $\mathfrak{g}(e)$ & $\mathfrak{g}_{0,ss}$ & $m_c$ & $\#\{dim(Z_{m_j})=1$\\
\hline
 $\mathfrak{su}(n,n)$ &$\mathfrak{sl}(2,\C)$ & $\mathfrak{sl}(n,\C)$& 1&0 \\
 \hline
$\mathfrak{so}(2,q)$ &$\mathfrak{sl}(2,\C)$ &$\mathfrak{so}(1, q-1)$ & 1&0\\
 \hline
 $\mathfrak{sp}(2n,\R)$ &$\mathfrak{sl}(2,\C)$&$\mathfrak{sl}(n,\R)$  & 1 &0\\
  \hline
 $\mathfrak{so}^*(4n)$ &$\mathfrak{sl}(2,\C)$&$\mathfrak{su}^*(2n)$ & 1 &0\\
 \hline
$\mathfrak{e}_7^{-25}$&$\mathfrak{sl}(2,\C)$ &$\mathfrak{e}_6^{-26}$ &  1&0 \\
 \hline
 $\mathfrak{f}_4^4$& $\mathfrak{g}_2$ &$\mathfrak{sl}(3,\R)$& 3 &2\\
 \hline
$\mathfrak{e}_6^2$ &$\mathfrak{g}_2$ & $\mathfrak{sl}(3,\C)$& 3&2\\
 \hline
 $\mathfrak{e}_7^{-5}$&$\mathfrak{g}_2$ &$\mathfrak{sl}(3,\mathbb{H})$ &3 &2\\
 \hline
$\mathfrak{e}_8^{-24}$ &$\mathfrak{g}_2$ & $\mathfrak{e}_6^{-26}$& 3&2\\
 \hline
 $\mathfrak{so}(p,q)$&$\mathfrak{so}(p-1,p)$&$\mathfrak{so}(1,q-p+1)$& $p-1$ &$p-1$\\
$2<p<q$&&&&\\
 \hline
 split & $\mathfrak{g}$ &$0$ & -&$rank(\mathfrak{g})$ \\
 \hline
 \end{tabular}
 \medskip
 \caption{Lie theoretic data for magical real forms, with the Lie algebras of the Cayley partner groups ($\mathfrak{g'}$) shown for the real forms of Hermitian tube type. }\label{Magicaldata}
 \end{table}

We can identify two extremes (see  Table  \ref{Magicaldata}) corresponding to $\mathrm{rank}(\mathfrak{g}(e))=1$ and  $\mathfrak{g}_{0,ss}=0$. These two extremes correspond to cases in which the Cayley component comes entirely from the first factor in \eqref{Cayleyproduct} (or more precisely from a factor of the form $\mathcal{M}_{K^{m_c+1}}(G')$, as described below) and, at the other extreme, to cases in which it comes entirely from the spaces of holomorphic differentials.

The first case, i.e.\ $\mathrm{rank}(\mathfrak{g}(e))=1$, corresponds to the real forms of Hermitian tube type.  In these cases $\mathfrak{g}(e)=\mathfrak{sl}(2,\C)$, $m_c=1$, and there is just one summand in $\bigoplus_{j=1}^{\mathrm{rank}(\mathfrak{g}(e))}H^0(K^{l_j+1}_{\Sigma})$. The Cayley component has the form 

\begin{equation}\label{Hsscayley}
\mathcal{M}_{K^{2}}(G_{0,ss}^\R)\times H^0(K^{2}_{\Sigma})
\end{equation}

Using the fact that $H^0(K^{2}_{\Sigma})$ can be identified with $\mathcal{M}_{K^{2}}(\R^+)$ we see that the Cayley components given by \eqref{Hsscayley} can be identified as moduli spaces of the form $\mathcal{M}_{K^{2}}(G')$, i.e.\ of $K^2$-twisted $G'$ Higgs bundles where $G'=G^{\R}_{0,ss}\times\R^+$.  The groups $G'$ are precisely the Cayley partner groups identified in \cite{BGRmaximalToledo}\footnotemark\footnotetext{The identification of the Lie algebras is clear. In order to relate the Cayley partners to the groups $G^{\R}_{0,ss}\times\R^+$ one must take care to consider groups with the same centers, for example $U^*(2n)=(\SU^*(2n)\ltimes\Z_2)\times\R^+$.} and we thus recover the Cayley correspondence for maximal components in $\Mhiggs G)$, as described in Section \ref{maxcomponents}.

The second case, i.e.\ $\mathfrak{g}_{0,ss}=0$, corresponds to the split real forms. In this case the first factor in the Cayley component does not vanish entirely but is just a finite set and we recover the Hitchin components described in Section \ref{HitchinFibComponents} as the images of Hitchin sections.

The other cases in Table \ref{Magicaldata} are a hybrid of the two extreme cases. The magical components in the associated moduli spaces are products of Cayley-like moduli spaces and Hitchin sections.  For real forms with the exceptional Lie algebras in the table, the Cayley components are all of the form

\begin{equation}\label{ECayley}
\mathcal{M}_{\Ksigma^4}(G_{0,ss}^\R)\times H^0(\Ksigma^2)\oplus H^0(\Ksigma^6)
\end{equation}

\noi where the Lie algebras of the groups $G_{0,ss}^\R$ are those shown in the table.

In the next section we describe the case of $\SO(p,q)$ with $2<p<q$ in some detail because it is illustrative, and also because the details play a role in Section \ref{global}. For full details see \cite{so(pq)BCGGO}.

\subsection{$\SO(p,q)$-Higgs bundles}\label{sopq}\hfill

\noi The magical $\mathfrak{sl}_2$-triples for which the Cayley real form is $\mathfrak{so}(p,q)$ (with If $2<p<q$) are principal in $\mathfrak{so}(2p-1,\C)\subset \mathfrak{so}(p+q,\C)$. They have
\begin{itemize}
\item $\mathfrak{g}(e)\simeq\mathfrak{so}(2p-1,\C)$,
\item $\mathfrak{g}_{0,ss}\simeq\mathfrak{so}(1, q-p+1)$, and
\item $C=S(O(1,\C)\times O(1,q-p+1))$ in $\SO(p,q)$
\end{itemize}

\noi The integers $\{m_j\}$ in the decomposition \eqref{sl2reps} are $\{1,3,\dots,2p-3\}\cup\{p-1\}$ with $m_c=p-1$.
This lead to Cayley maps of the form 

 \begin{equation}
        \label{EQ PsiOnModuli}
        \Psi:\cM_{K^p}(\Sigma, \SO(1,q-p+1))\times \displaystyle\bigoplus\limits_{j=1}^{p-1}H^0(K^{2j})\longrightarrow\cM(\Sigma, \SO(p,q)),
     \end{equation} 

We can explicitly describe these maps and the structure of the resulting Cayley components. We start with a convenient description of  $\SO(p,q)$-Higgs bundles.  The maximal compact subgroup of $\SO(p,q)$ is $H=S(O(p)\times O(q)$.  According to Definition \ref{GHiggsDefn} the bundle in an $\mathrm{SO}(p,q)$-Higgs bundle should thus be a principal $S(O(p,\C)\times O(q,\C)$-bundle. We may use the standard representations of $O(p,\C)$ and $O(q,\C)$ to give a vector bundle description as follows.

\begin{definition}\label{DEF SO(p,q) Higgs bundles}
An {\bf $\mathrm{SO}(p,q)$-Higgs bundle} on a closed surface $\Sigma$ is a tuple $((V,Q_V),(W,Q_W),\eta)$, where $(V,Q_V),$ $(W,Q_W)$ are rank $p$ and rank  $q$ holomorphic orthogonal vector bundles respectively such that $\det(V\oplus W)=\mathcal{O}$,  and $\eta\in H^0(\Hom(W,V)\otimes \Ksigma)$.
\end{definition}

\begin{remark}
Unless it is essential we often suppress the orthogonal structures and denote the $\mathrm{SO}(p,q)$-Higgs bundles by tuples $(V,W,\eta)$.
\end{remark}

Definition \ref  {DEF SO(p,q) Higgs bundles} applies also if $p=1$ and the Higgs fields are $L$-twisted, where $L$ is any line bundle. In particular, a {\bf $\Ksigma^m$-twisted $\mathrm{SO}(1,n)$-Higgs bundle} is a tuple $(\mathcal{I}, W,\eta)$ where 

\begin{itemize}
\item $\mathcal{I}^2=\mathcal{O}$, i.e.\ $\mathcal{I}$ denotes a square root of the trivial line bundle,
\item $W$ is an orthogonal bundle of rank $n$,
\item $\eta$ is a holomorphic map $\eta: W\rightarrow \mathcal{I}\otimes \Ksigma^{m}$
\end{itemize}

The Cayley map \eqref{EQ PsiOnModuli} can then be explicitly described (see Section 4 in \cite{so(pq)BCGGO}) by

\begin{equation}
\Psi: \big ([\mathcal{I},\widehat W,\eta], q_2,q_4,\dots q_{2p-2}\big)\longmapsto [V,W,\eta]
\end{equation}

\noi with 
\begin{equation} \label{CayleyHitchin}
(V,W,\eta)= \left(\mathcal{I}\otimes \mathcal{K}_{p-1},\ \mathcal{I}\otimes \mathcal{K}_{p-2}\oplus \widehat W,\  \begin{bmatrix}q_2&q_4&\cdots& q_{2p-2}&\eta_{\widehat W}\\1&q_2&\cdots&q_{2p-4}&0\\
\vdots &\ddots&\ddots&\vdots&\vdots\\
 0&\cdots&1&q_2&0\\0&\cdots&0&1&0\end{bmatrix}\right)
\end{equation}

\noi where
\begin{itemize}
\item $(\mathcal{I}^2,\widehat W,\eta_{\widehat W})$ defines a $\Ksigma^p$-twisted $\SO(1,q-p+1)$ Higgs bunde,
\item $q_j\in H^0(\Ksigma^{j})=H^0(Hom(\Ksigma^{i},\Ksigma^{i+j-1})\otimes \Ksigma)$ for $j=2,4,\cdots,2p-4$, and
\item $\mathcal{K}_n=\Ksigma^{n}\oplus \Ksigma^{n-2}\oplus\cdots\oplus \Ksigma^{2-n}\oplus \Ksigma^{-n}$.
\end{itemize}

In  \cite{so(pq)BCGGO} we showed directly that the $\SO(p,q)$-Higgs bundles of the form in \eqref{CayleyHitchin} comprise a union of connected components in $\Mhiggs G)$. In order to relate this to the general results for Cayley components coming from magical $\mathfrak{sl}_2$-triples we proceed as follows.

Let $\SO(p,q)\subset \SO(p+q,\C)$ be the group of unit determinant block matrices $\begin{bmatrix}A&B\\C&D\end{bmatrix}$ such that

\begin{equation}
\begin{bmatrix}A^T&C^T\\B^T&D^T\end{bmatrix}\begin{bmatrix}Q_p&0\\0&-Q_q\end{bmatrix}\begin{bmatrix}A&B\\C&D\end{bmatrix}=\begin{bmatrix}Q_p&0\\0&-Q_q\end{bmatrix}
\end{equation}

\noi where $Q_p$ and $Q_q$ are positive definite symmetric bilinear forms on $\C^p$ and $\C^q$ respectively

\noi Without loss of generality we assume that $p\le q$ Then we can identify the subgroup $C=S(\mathrm{O}(1,\C)\times\mathrm{O}(q-p+1,\C))$ embedded via

\begin{equation}
(\det(A),A)\mapsto\begin{bmatrix}\det(A)I_p&0&0\\ 0&\det(A)I_{p-1}&0\\0&0&A\end{bmatrix}
\end{equation}

\noi and a copy of $T=\C^*$ embedded via

\begin{equation}
t \mapsto\begin{bmatrix}T_p(t)&0&0\\0&T_{p-1}(t)&0\\0&0&I_{q-p+1}\end{bmatrix}
\end{equation}

\noi where
\begin{equation}
T_n(t)= \begin{bmatrix}t^{n-1}&0&\cdots&0\\
0&t^{n-2}&\cdots&0\\
\vdots&\ddots&&0\\
0&\cdots&\cdots&t^{1-n}\\
\end{bmatrix}
\end{equation}

\noi Let $\mathcal{E}_C$ be the frame bundle for $\mathcal{I}\oplus\widehat{W}$, and let $\mathcal{E}_T$ be the frame bundle for $\Ksigma$.  We can then write

\begin{equation}
\big (\mathcal{I}\otimes \mathcal{K}_{p-1}\big ) \oplus\Big ( \big (\mathcal{I}\otimes \mathcal{K}_{p-2}\big )\oplus\widehat W\Big )\simeq\mathcal{E}_C\star\mathcal{E}_T[H^{\C}][\C^p\oplus\C^q]
\end{equation}

\noi where the $\star$ denotes the operation coming from the product of the subgroups $C$ and $T$ in $S(O(p,\C)\times O(q,\C))\subset\SO(p+q,\C)$.  This shows that the frame bundle for vector bundle $V\oplus W$ is precisely of the special form required in \eqref{specialHiggs}.   We omit the details showing how to realize the Higgs fields in \eqref{CayleyHitchin} as Higgs fields of the form in \eqref{specialHiggs} except to note that 
\begin{itemize}
\item the contribution from 1's corresponds to $f$, 
\item for $i=2,4,\dots,2p-2$  the contribution from the differentials $q_i$ corresponds to $\phi_{m_{i-1}}$, and
\item the contribution from $\eta_{\hat W}$ corresponds to $\phi_{m_{p-1}}$.
\end{itemize}

\subsection{Local minima of the Hitchin function}\hfill

\noi   We note that  (see  Theorem 5.10 in \cite{so(pq)BCGGO})

\begin{theorem}In \eqref {CayleyHitchin} if the holomorphic differentials $q_i$ are all zero and $[\mathcal{I},\widehat W,\eta]$ represents a local minimum of the Hitchin function \eqref{fdef} on $\Mhiggs \SO(1,q-p+1))$ then the resulting $\SO(p,q)$-Higgs bundles define local minima of the Hitchin function \eqref{fdef}. These are the only local minima with non-trivial  Higgs fields. 
\end{theorem}

If $q-p+1>2$ then the local minima on $\Mhiggs \SO(1,q-p+1))$ are all trivial, i.e.\ have $\eta_{\widehat W}=0$. In contrast, $\Mhiggs \SO(1,2))$ admits non-trivial minima (see Section \ref{sopqextra})

All possible topological types of $H^{\C}$-bundles occur in $\Mhiggs \SO(p,q))$ and trivial local minima occur with bundles of all such topological types.  This, together with a count for $|\pi_0(\Mhiggs \SO(1,n))|$, leads to the count shown in Tables \ref{SO2nTable} and \ref{SO2n+1Table} of connected components in $\Mhiggs \SO(p,q))$. 

The factors $\prod_{j=1}^{\mathrm{rank}(\mathfrak{g}(e))}H^0(K^{l_j+1}_{\Sigma})$ in \eqref{EQ PsiOnModuli} are contractible so the precise number of components in the image of the Cayley map (called  {\it Cayley components}) is determined by $\pi_0(\mathcal{M}_{K^{m_c+1}}(G_{0,ss}))$. This in turn depends not only on Lie algebra data but also on the center of $G_{0,ss}$.  
The tables list cases where the component count has been done.



\section{Remaining cases and Loose ends} 

There are two types of real forms that have not yet been discussed: the non-tube Hermitian forms and those in the last column of Table \ref{realforms}.  The non-tube Hermitian forms can be understood in terms of their tube-type cousins - see Section \ref{non-tube}.  For the real forms in the last column, the classical cases have been analyzed individually. There are no published results for the exceptional cases in the last column but it is reasonable to conjecture that for all of them the only local minima of the Hitchin function are at $\Phi=0$.   

In addition to the above cases, we discuss the other remaining situations in which the tally of connected components is (potentially) incomplete. 

\subsection{Non-tube type}\label{non-tube}

The Lie algebras $\mathfrak{su}(p,q)$ with $p\ne q$, $\mathfrak{so}^*(4m+2)$, and $\mathfrak{e}_6^{-14}$ are Hermitian real forms of non-tube type.   They are not on the list in Theorem \ref{classification}, but the results of section \ref{magicalSL2} apply in modified form. In each case there is a {\it maximal tube subalgebra}, denoted $\mathfrak{g}_T\subset\mathfrak{g}$, namely $\mathfrak{su}(p,p), \mathfrak{so}^*(4m)$, and $\mathfrak{so}(2,8)$ (see Section 6 in \cite{BGRmaximalToledo}).  

If $G_T\subset G$ is the subgroup with Lie algebra $\mathfrak{g}_T$, and $G_T^{Ad}$ is its adjoint form, then (see Theorem 6.2. \cite{BGRmaximalToledo}) the maximal components in the moduli spaces for $G$ and $G_T$ are related by a fibration

\begin{equation}\label{maxtubefibration}
M (L_{\C})\rightarrow \Mhiggs G)_{max} \rightarrow \Mhiggs G_T^{Ad})_{max} 
\end{equation}

\noi where $M(L_{\C})$ is the moduli space of polystable $L_{\C}$-bundles. If $G=\SU(p,q)$ then $L_{\C}=\GL(q-p,\C)\rtimes\Z_p$  and $L_{\C}=\C^*$ in the other cases.

One consequence of \eqref{maxtubefibration} is that $\Mhiggs G)_{max}$ has lower than expected dimension. Moreover, since $M(L_{\C})$ is connected it also shows that $\Mhiggs G)_{max}$ has the same number of components as $\Mhiggs G_T^{Ad})_{max} $. For $G=\SU(p,q)$ or $\SO^*(4m+2)$ this does not yield new information. In the case of $E_6^{-14}$, where $L=U(1)$ and $G_T=Spin_0(2,8)$ this leads to the result that the $\Mhiggs E_6^{-14})_{max}$  has the same number of connected components as $\Mhiggs \PSO_0(2, 8))_{max}$. 
As shown in \cite{GP-2017} (Theorem 5.6), this number is 2.

\subsection{Non-maximal components}\label{nonmax}

For Hermitian real forms the Cayley components are maximal, i.e.\ the Toledo invariant is maximal on these components. A case by case analysis has shown that these are the only maximal components. It is also known that there is at least one non-empty connected component for every other value of the Toledo invariant strictly within its bounded range.  The count depends on the number of connected components in the locus of local minima for the Hitchin function on these components.  If $G=\SU(p,q)$ or $\PU(p,q)$ these non-maximal components are known to be connected, and the non-maximal components have been counted in the cases of $G=\SO(2,2),\ \SO(2,3),\ \SO_0(2,3)=\PSp(4,\R)$, and $ \Sp(4,\R)$ (see the references in the tables).  For the other classical matrix groups of Hermitian type the connectivity of the non-maximal components remains an open question.

\subsection{Non-magical real forms}  The real forms in the last column of Table \ref{realforms} do not admit magical $\mathfrak{sl}$-triples and are not Hermitian of non-tube type. 

The case $\SO(1,n)$ was analyzed in \cite{AG11} where it was shown that the Hitchin function has only trivial local minima of on $\Mhiggs \SO(1,n))$ and hence that the moduli space is connected. In \cite{Oliveira_GarciaPrada_2016} they show that the Hitchin function has only trivial local minima of on $\Mhiggs U^*(2n))$. Since the maximal compact subgroup of $U^*(2n)$, i.e.\ $\mathrm{Sp}(2n)$, is connected, it follows that $\Mhiggs U^*(2n))$ is connected. The maximal compact subgroup for $U^*(2n)$ is the same as for $\SU^*(2n)$ and hence it follows that $\Mhiggs \SU^*(2n))$ is also connected. 

There is reason to conjecture that if $G$ is a real Lie group whose Lie algebra is one of the exceptional real forms in this category, namely $\mathfrak{e}_6^{-26}$ and $\mathfrak{f}_4^{-20}$,  then all local minima on $\mathcal{M}(\Sigma,G)$ lie in $f^{-1}(0)$, i.e.\ have $\Phi=0$ (see \cite{collier2023}).

Inspection of Table \ref{realforms} shows that for all other cases there are local minima with $\Phi\ne 0$.  For the classical matrix groups the case-by-case analyses show that the the non-trivial local minima (and hence the components on which they occur) are either Cayley components or components with non-maximal Toledo invariant. 
The latter occur only in the moduli spaces of $G$-Higgs bundles for $G$ of Hermitian type.  

\subsection{Non-magical companions to Cayley components}

The Cayley map detects connected components that are not accounted for simply by the topological types of the principal bundles. Via the NAH correspondence these components correspond to connected components in $\mathrm{Rep}(\pi_1(S),G)$ in which the representations cannot be deformed to representations which factor through a maximal compact subgroup of $G$. These features follow from the fact that the Higgs bundles in these components cannot be deformed to $G$-Higgs bundles with zero Higgs field. 

This leaves open the possibility that in addition to the Cayley components there are other unaccounted for connected components in which the Higgs fields are never identically zero.  This is known not to be the case if $G$ is a split real form of a classical group or if $G=\SO(p,q)$. For the real forms of hemitian tube-type it is known that there are non-maximal components but, as discussed in Section \ref{nonmax}, the precise number is not known in all cases.  

For the exceptional real forms in Table  \ref{realforms}  the issue could be settled if it were known that the Hitchin function has no local minima other than at zero or those which occur in the Cayley components.   The analysis of local minima can be separated into two cases, namely smooth local minima and local minima outside the smooth locus of the moduli space. On the smooth locus the criterion in theorem \ref{mincrit} can be used to identify necessary conditions that must be satisfied at a smooth minimum.  Showing that these cannot be satisfied except at known minima would rule out the possibility of undiscovered smooth local minima.   The possibility of local minima at singular points has to be ruled out differently.


\section{Beyond $\pi_0$}

The local minima of the Hitchin function identify connected components on $\Mhiggs G)$, and the Cayley map reveals partial structural information about certain distinguished components. In some cases, which we describe below, the global structure can be described much more explicitly. The descriptions lead to global information which does not depend on the complex structure on $\Sigma$ and thus applies to the representation varieties $\mathrm{Rep}(\pi_1(S),G)$ on the other side of the NAH correspondence.

\subsection{Global descriptions}\label{global}  In \cite{selfduality} Hitchin described all components of the moduli space $\Mhiggs \SL(2,\R))$. Similar methods were used in \cite{bradlow-garcia-prada-gothen:2012,gothen-01} to describe the maximal components in $\Mhiggs \mathrm{Sp}(4,\R))$, and in \cite{gothen-01} to describe components in $\Mhiggs P\mathrm{Sp}(4,\R))$.  All of these examples are essentially special cases of the description in \cite{CollierSOnn+1components}
 of distinguished components of $\Mhiggs \SO(n,n+1))$, as can be seen because of the isomorphisms $\PSL(2,\R)\simeq\SO_0(1,2)$ and $P\mathrm{Sp}(4,\R)\simeq\SO_0(2,3)$.  
 
Though not obtained in this way, the special components described in \cite{CollierSOnn+1components} are in fact Cayley components. The $\SO(n,n+1)$-Higgs bundles in these components are all of the form in \eqref{CayleyHitchin}, but since $q-p+1=2$, the orthogonal bundles $\widehat W$ in the $\Ksigma^p$-twisted $\SO(1,q-p+1)$ Higgs bundles are $O(2,\C)$-bundles. Such bundles are classified topologically by first and second Stiefel-Whitney classes $(sw_1,sw_2)\in H^1(\Sigma,\Z_2)\times H^2(\Sigma,\Z_2)$.  If $sw_1=0$ then $\widehat W$ is an $\SO(2,\C)$-bundle which may be taken to be of the form $L\oplus L^{-1}$, where $L$ is a line bundle of non-negative degree. The degree, $d$, of $L$ is bounded by $n(2g-2)$ and Collier shows in \cite{CollierSOnn+1components} that

\begin{theorem} [Theorem 4.11 in \cite{CollierSOnn+1components} ] For each degree in the range $0<d\le n(2g-2)$ the resulting component is diffeomorphic to a product

\begin{equation}\label{dcomponent}
\mathcal{F}_d\otimes\bigoplus_{j=1}^{n-1}H^0(\Ksigma^{2j})
\end{equation}

\noi where $\mathcal{F}_d$ is the total space of a rank $d+(2n-1)(g-1)$ vector bundle over the symmetric product $Sym^{n(2g-2)-d}(\Sigma)$.  In particular
\end{theorem}

\begin{corollary} The components given by \eqref{dcomponent} are smooth and deformation retract onto the symmetric product $Sym^{n(2g-2)-d}(\Sigma)$. 
\end{corollary}

\begin{theorem}If $d=0$ then the resulting component is homeomorphic to a product

\begin{equation}
\mathcal{F}_0\otimes\bigoplus_{j=1}^{n-1}H^0(\Ksigma^{2j})
\end{equation}

\noi where $\mathcal{F}_0$ is a singular space which may be described as a GIT quotient $\tilde{\mathcal{F}}_0//O(2,\C)$ of the space

\begin{equation}
\tilde{\mathcal{F}}_0=\big\{(M,\mu,\nu)\ |\ M\in Pic^0(\Sigma),\ \mu\in H^0(M^{-1}\Ksigma^n), \ \nu\in H^0(M\Ksigma^n)\big\}
\end{equation}
\end{theorem}
 
The components in which the first Stiefel-Whitney class of $\widehat W$ does not vanish can also be described (see Theorem 5.3 in \cite{CollierSOnn+1components}) and are of the form

\begin{equation}\label{sw12components}
\mathcal{F}_{sw_1}^{sw_2}\otimes\bigoplus_{j=1}^{n-1}H^0(\Ksigma^{2j})
\end{equation}

\noi with $\mathcal{F}_{sw_1}^{sw_2}=\tilde{\mathcal{F}}_{sw_1}^{sw_2}/(\Z_2\oplus\Z_2)$ where $\tilde{\mathcal{F}}_{sw_1}^{sw_2}$ is a rank $(4n-2)(2g-2)$ vector bundle.  The base of this vector bundle is one of the connected components of the Prym variety for the double cover of $\Sigma$ determined by $sw_1$, with the component determined by $sw_2$. 

\begin{theorem} The components given by \eqref{sw12components} deformation retract onto the quotient of the Prym variety by the $\Z_2$-action coming from the involution in the double cover.
 \end{theorem}
 
 \subsection{Other topological information}  In \cite{bradlow-garcia-prada-gothen:2008} we exploit Morse-theoretic information at higher critical points of the Hitchin function to compute higher homotopy groups (i.e.\ $\pi_i$ with $i>0$) for components of the Higgs bundle moduli spaces $\Mhiggs \GL(n,\C))$ and $\Mhiggs U(p,q))$.  In \cite{AndrePGLnR} the author proves, also by Morse theoretic methods,  that the space consisting of two of the components in $\Mhiggs \SL(3,\R))$ is homotopically equivalent to the moduli space of polystable $O(3,\C)$ bundles. The third component is the contractible Hitchin component.


 \section{Supplemental Notes}

The tables in Section \ref{Tables} summarize known information about the number of connected components in $\Mhiggs G)$ if $G$ is a real form with Lie algebra in Table \ref{realforms}. 

The number of connected components in $\Mhiggs G)$ is determined not just by the Lie algebra of $G$ but also by the center of the group. This explains why the entries in the tables differ for groups with the same Lie algebra but different centers. 

The entries in the table are obtained, for the most part, from the techniques outlined in this survey. The references contain the details but we record here some additional explanatory notes to give a flavor of the arguments or to explain some stand-out cases.  The entries in [[brackets]] have been computed as lower bounds but are conjectured to be precise. The uncertainty comes from gaps in our knowledge concerning the non-trivial local minima of the Hitchin function, especially in components with non-maximal Toledo invariant, or outside the smooth locus of the components. 
It is reasonable to conjecture that the only connected components not labeled by topological invariants of the principal bundles, i.e.\ with non-trivial local minima of the Hitchin function, are Cayley components. 

We organize this section differently to the way the tables are organized to emphasize the commonalities in the features which determine component counts in $\Mhiggs G)$ for the various categories of real forms.

\subsection{Real forms of Hermitian symmetric type}\hfill

\medskip 
\subsubsection{$\bm{\SL(2,\R)(\simeq\Sp(4,\R)\simeq\SU(1,1))}$:} The $\SO(2)$-representations (i.e.\ where $\Phi = 0$) correspond to stable {\it degree zero} line bundles, i.e.\ correspond to {\it zero Euler class} (aka Toledo invariant).  The automorphism which switches $L\oplus L^{-1}\mapsto L^{-1}\oplus L$ lies in $O(2)$ but not in $SO(2)$. This explains why we get one component for each integer in $[1-g,g-1]$.

\medskip 
\subsubsection{$\bm{\SO^*(2n)\  \mathrm{and}\ \mathrm{Sp}(2n,\R)}$:} The lower bound includes the components at $\pm\tau_{max}$ and reflects the fact that

\begin{enumerate} 
\item there are no `topological' components if the relevant Toledo invariant lies between zero and its extremal values. [The Higgs bundles are defined by triples $(V,\beta,\gamma)$ where $V$ is a $U(n)$-bundle of degree $d$; if $d\ne 0$ then $V$ does not support flat connections], and
\item there is at least one component for each non-zero, non-extremal value of  the Toledo invariant (detected by the non-trivial minima of the Morse function), 
\item For $\SO^*(2n)$ there is just one component with zero or maximal Toledo invariant ($\lfloor\frac{n}{2}\rfloor(2g-2)$). The total number of components is thus bounded below by the number of values for the Toledo invariant, with equality if and only if each of the non-maximal components is connected.
\end{enumerate}

Thus the lower bound is exact unless there are undetected exotic components with non-zero non-extremal Toledo invariant. Note: the $\SO^*(4n)$-moduli spaces have  magical components but the $\SO^*(4n+2)$ do not.

\medskip 
\subsubsection{$\bm{\PSO^*(2n)}$}  The maximal compact subgroup of $\PSO^*(2n)$ is $U(n)/\Z_2$. Topological types of principal $U(n)/\Z_2$-bundles are classified by $(d,w)\in\Z\times \Z_2$ if $n$ is even and $d\in\Z$ if $n$ is odd. Not all $U(n)/\Z_2$-bundles lift to $U(n)$.

For $n=2m$, it can be shown (see \cite{GP-2017}) that $0\le |d|\le 2m(g-1)$ and that $U(n)/\Z_2$-bundles lift to $U(n)$ iff $w=0$. Using the Cayley correspondence one can show that there are two maximal components, one which lifts to $\SO^*(2m)$, and one which does not.  There are two components in which $\Phi=0$, namely those labelled by $(d,w)=(0,w)$.

For $n=2m+1$ the invariant $d$ must satisfy $0\le |d|\le 2m(g-1)$. The $U(n)/\Z_2$-bundles lift to $U(n)$ iff $d$ is even, in particular, if $|d|$ is maximal. There is no obstruction to lifting a principal bundle from $\PSO^*(2m+1)$ to $\SO^*(2n)$ (See also \cite{GW10}). There is thus one maximal component since there is just one maximal component in $\Mhiggs, \SO^*(2m+1)$. There is one component in which $\Phi=0$, namely the one labelled by $d=0$.

Note:  $\PSO^*(2)$ is the trivial group, while  $\PSO^*(4)$ and $\PSO_0(2, 2)$ are not simple.

\medskip 
\subsubsection{ $\bm{\SU(p,p)}$}  The components detected by the non-trivial minima are all connected, except for the maximal $p=q$ case in which we get $2^{2g}$ components (because the fixed determinant condition is on $det(E)^2$). The number of non-maximal values taken by the Toledo invariant is

\begin{equation}
2(\tau_{max}-1)+1=2(p(g-1)-1)+1=2p(g-1)-1
\end{equation}

\noi  Moreover if $p<q$ then the maximal components fiber over $M_0(U(q-p))$ with fiber $\Mhiggs \SU(p,p))_{\tau_{max}}$, so the component count is the same as for $\SU(p,p)$.  Here $M_0(U(q-p))$ is the moduli space of degree zero, rank $q-p)$ polystable vector bundles.

Note that there are no `topological' components if the Toledo invariant is non-zero (since vector bundles with non-zero degree  cannot support flat connections). 

In \cite{gothen-01} Gothen shows that for $\SU(2,2)$ the components with Toledo invariant equal to zero or maximal are connected.

\medskip 
\subsubsection{$\bm{\mathrm{Sp}(4,\R)}$:} The bound on $d$ is $|d|\le 2g-2$. There is one component ($d=0$) where $\Phi=0$, two each of the multiple components with $d=\pm (2g-2)$, and one component for each $0<|d|<2g-2$. The maximal components include the so called Gothen components (first identified in \cite{gothen-01}) and the Teichm\"uller  components.

\medskip 
\subsubsection{$\bm{\mathrm{Sp}(2n,\R)}$:} The bound on $d$ is $|d|\le n(g-1)$.  There is one component ($d=0$) where $\Phi=0$. There are two of each maximal component ($d=\pm n(g-2)$) and for each maximal value there are $2\cdot 2^{2g}$ Cayley components with trivial minima and $2^{2g}$ Teichm\"uller  components. There is at least one component for each $0<|d|<n(g-1)$.

\medskip 
\subsubsection{$\bm{P\mathrm{Sp}(2n,\R)}$:} See \cite{GP-2017}. 
\begin{itemize}
\item $n$ even: Topological types of principal bundle are determined by $(d,w)\in \Z\times\Z_2$ with $|d|\le n(g-1)$, for a total of $2\cdot \big (2n(g-1)+1\big )=4n(g-1)+2$  types. Two types have $d=0$ and two have $\pm$ maximal value. Each maximal value accounts for $2^{2g+1}+2$ components (if $n\ge 4$). $2^{2g}$ of the $\PSp(2n,\R)$-components do not lift to $\mathrm{Sp}(2n,\R)$.
\item $n$ odd: There is no obstruction to lifting from $P\mathrm{Sp}(2n,\R)$ to  $\mathrm{Sp}(2n,\R)$ (See also \cite{GW10}).
Topological types of principal bundle are determined by $d\in \Z$ with $|d|\le 2n(g-1)$, for a total of $4n(g-1)+1$  types. One type has $d=0$ and two have $\pm$ maximal value. Each maximal value accounts for $3$ components (if $n\ge 3$).  The $3.2^{2g}$ maximal $\Sp(2n,\R)$-components map $2^{2g}:1$ onto the maximal $\PSp(2n,\R)$-components
\end{itemize}

\medskip 
\subsubsection{$\bm{\SO(2,q)}$}\label{SO2q}  The maximal compact subgroup of $\SO(2,q)$ is $H=S(O(2)\times O(q))$. Thus the primary topological invariants are given by  $(a,b,c)\in H^1(S,\Z_2)\times H^2(S,\Z_2)\times H^2(S,\Z_2)$.

If  the invariant $a\in H^1(S,\Z_2)$ is zero then structure group of the Higgs bundles reduces to $\SO_0(2,q)$. In this case the group is of Hermitian type and a new invariant appears (in this case, the new invariant comes from a lift of $b\in H^2(S,\Z_2)$ to $H^2(S,\Z)\simeq\Z$).   The primary topological invariants are then $(l,c)\in H^2(S,\Z)\times \times H^2(S,\Z_2)$ with $|l|\le 2g-2$.  The integer $l$ is, up to a normalizing factor, the Toledo invariant mentioned in Section \ref{maxcomponents}. However  in the group $SO(2,q)$ the sign of the Toledo invariant is not well defined, or said another way, positive and negative Toledo invariant for $\SO_0(2,q)$ can be conjugated to each other using $SO(2,q)$. When $|l|$ is maximal new invariants appear, with $c$ being replaced by a pair in $\Z^{2g}\times\Z_2$.

\begin{itemize}
\item If $a\ne 0$ then all components admit only trivial minima, giving $(2^{2g}-1)\cdot 2\cdot 2$ such components.
\item If $a=0$, so that the invariants become $(|l|,c)$ as above then there are
\begin{itemize}
\item  two components with only trivial minima (with $l=0$), 
\item $2^{2g}$ Cayley components (with $|l|=2g-2$), and 
\item at least $2\cdot (2g-3)$ other components (with $0<|l|<2g-2$).
\end{itemize}
\end{itemize}

\medskip 
For $\bm{\SO(2,2)}$ if $a\ne 0$ then the other two invariants remain $\Z_2$-valued.  The non-zero minima in $\Mhiggs \SO(2,2)$ all lie in components with $a=0$, i.e.\ where we have $\SO_0(2,2)$-Higgs bundles.  In that case the maximal compact subgroup of $\SO_0(2,2)$ is $H=SO(2)\times SO(2)$ and the invariant $c\in H^2(S,\Z_2)$ also lifts to $H^2(S,\Z)$.  The primary topological invariants are  $(l,m)\in H^2(S,\Z)\times \times H^2(S,\Z)\simeq \Z\times\Z$, but subject to a Milnor-Wood bound  which in this case is $l\ge 0$ and 

\begin{equation} l-2g+2\le m\le -l+2g-2\ .\end{equation}

\noi Also,  if $l=0$, then only $|m|$ is an invariant.

\medskip 
For $\bm{\SO(2,3)}$ the extra Cayley components come from the fact that  if $q=3$ then the $\SO(q)$-bundle reduces to an $S(O(1)\times O(2))$-bundle. Moreover, if the structure group reduces further to $\SO(1)\times \SO(2)$ then the second Stiefel-Whitney class is replaced by a new integer invariant.

\medskip 
For $\bm{\SO_0(2,2)}$ all the minima in $\Mhiggs g)$ are connected subvarieties (see Section 6 of \cite{so(pq)BCGGO}) except when $(l,m)$ equals $(0,2g-2)$ or $(2g-2,0)$; in each of those case there are $2^{2g}$ Hitchin components.  There is one component (when $(l,m)=(0,0)$) with $\Phi=0$ at a local minimum.  

Note:  The $2^{2g}$ choices in the Hitchin components corresponds to an allowed twisting by a square root of $\mathcal{O}$ in both of the $\SO(2,\C)$ bundles in the Higgs bundle. The constraint on the degree zero bundle is required to preserve the trivial determinant condition.

\medskip 
For $\bm{\SO_0(2,3)=P\mathrm{Sp}(4,\R)}$, in \cite{AC} it is shown that there are  $2(2^{2g} -1) + 4g -3$ connected components for each of the maximal values of the integer topological invariant. In \cite{GO12} it is proven that there is one connected components for each other values of the invariants. This gives $2(2(2^{2g} -1)+4g-3)+4(2g-3)+2=2^{2g+2} +16g-20$ connected components. The extra maximal components (compared to the number for $\SO_0(2,q), q>3$) arise in the same way as the extra Cayley components for $\SO(2,3)$, except now the components labelled by $\pm d$ must be counted separately.  There are only two components in which the local minimum of the Hitchin function is at points with $\Phi=0$, namely those with  $a=l=0$ as above. 
\medskip 
\subsubsection{$\bm{\PSO_0(2,2n)}$}\hfill

The maximal compact subgroup is $(\SO(2)\times\SO(2n)/\Z_2$, leading to a classification of the principal bundles by $(d,w)\in\Z\times\Z_2$ with $0\le |d|\le 4g-4$ (see \cite{GP-2017}). The bundles with even values of $d$ lift to $\SO(2)\times\SO(2n)$-bundles. In particular, all bundles in the maximal components lift, so the maximal $\SO_0(2,2n)$-components surject onto $\Mhiggs \PSO(2,2n)$. There are $2^{2g}$ such $\SO_0(2,2n)$-components for each $w\in\Z_2$, but these get identified in $\Mhiggs \PSO(2,2n)$. The two components with $d=0$ are the only components where the Higgs field can vanish.

\subsection{Split real forms}\hfill

\medskip
 

\subsubsection{$\bm{\mathrm{\PSL}(2,\R)}$}  In \cite{goldman-88} there is a formula for the number of components of $Hom(\pi_1,G)$  where $G$ is any $n$-fold cover of $\PSL(2,\R)$, viz.:

\begin{align}
2n^{2g}+\frac{4g-4}{n}-1\ &\mathrm{if} \ n\ \mathrm{divides}\ 2g-2\\
2\Big[\frac{2g-2}{n}\Big]+1\ &\mathrm{otherwise}
\end{align}

\noi $\PSL(2,\R)$ can be identified with the identity component of $\PGL(2,\R)$. The latter has two components distinguished by the sign of the determinant so the two components of $\GL(2,\R)$ remain distinct in $\PSL(2,\R)$).

\medskip 
\subsubsection{$\bm{\mathrm{PSL}(n,\R)}$}  In \cite{liegroupsteichmuller} The description of components for $\mathrm{PSL}(n,\R)$  proceeds as follows:

\begin{itemize}
\item If $n$ is odd, $\mathrm{PSL}(n,\R)= \SL(n, \R)$ and if $w_2(E) = 0$ there are two components, if $w_2(E)\ne 0$ just one.

\item If $n = 2m$, there are four topological types. Two do not lift to $\SL(n, \R)$ and so give connected spaces. Two do lift and if $w_2(E) \ne (g - l)m^2$ this gives a connected space. If $w_2(E) = (g - l)m^2$ we have three components-the one containing minima with $\Phi=0$ and two copies 
 \footnotemark\footnotetext{The two copies get identified inside the $\PSL(2n,\C)$-moduli space by an isomorphism in the non-identity component of $\PGL(2n,\R)$ (see the last paragraph in \cite{liegroupsteichmuller}).} of the Teichm\"uller  component. 
\end{itemize}

\medskip
 

\subsubsection{$\bm{\mathrm{\PGL}(2n,\R)}$} Two factors complicate the counting of connected components of $\Mhiggs \PGL(2n,\R))$, namely (a) the maximal compact subgroup, $\mathrm{PO}(2n)$ has two components, and (b) $\mathrm{PO}(2n)$ is not a matrix group. The first factor affects the number of topological types of principal bundles, while the second means that vector bundle methods are not immediately available.  Both are dealt with in \cite{AndrePGLnR}.


\subsection{The real forms $\SO(p,q)$}\label{sopqextra}\hfill

\medskip 

The maximal compact subgroup of $\SO(p,q)$ is $H=S(O(p)\times O(q))$. Assuming $2<p\le q$,  the primary topological invariants are thus 

$$(a,b,c)\in H^1(S,\Z_2)\times H^2(S,\Z_2)\times H^2(S,\Z_2)\ $$

\noi The invariant $a$ distinguishes between $\SO(p,q)$- and $\SO_0(p,q)$-Higgs bundles, where $\SO_0(p,q)$-Higgs denotes the connected componet of the identity, with $a=0$ corresponding to the latter.

In Theorem 5.10 of \cite{so(pq)BCGGO} (see also theorem \ref{CayleyHitchin}) we classify all the local minima on $\mathcal{M}(SO(p,q))$ for all $1\le p\le q$. This includes one $\Phi=0$ minima on each of the $2^{2g+2}$ non-empty component of the moduli space of polystable $S(O(p,\C)\times O(q,\C))$-bundles. The other minima, i.e.\ those which detect the Cayley components, lie in sectors with topological invariants listed in the table in \S 6.1 of \cite{so(pq)BCGGO}. Except for $\SO(p,p+1)$, in all other cases this yields a total of $2^{2g+1}$ extra Cayley/exotic components.  

If $p$ is even then all Cayley components lie in sectors with $a=0$, while if $p$ is odd then just two of the Cayley components have $a=0$. These considerations account for the dependence on the parity of $p$ in the component counts for $\Mhiggs \SO_0(p,q))$.

\medskip 
\subsubsection{$\bm{\mathrm{\SO}(n,n+1)}$}  In the case of $\SO(p,q)=\SO(n,n+1)$ with $n\ge 3$ there are an additional $2n(g-1)-1$ so-called Collier components, first described in \cite{CollierSOnn+1components}, due to the fact that $q-p+1=2$ if $q=p+1$  (see Section \ref{global}).
\medskip 
\subsubsection{$\bm{\mathrm{\SO}(1,2)}$} For $\SO(1,2)$,  $2(2^{2g}-1)$ of the topological components come from the classification of $O(2,\C)$ bundles with $w_1\ne 0$, and one additional topological component comes from the component in which the bundles are of the form $(I,W)=(\mathcal{O}, L\oplus L^{-1})$ with $deg(L)=0$. The other components have $0<deg(L)\le 2g-2$ and non-trivial Higgs field.
\medskip 
\subsubsection{ $\bm{\SO(1,1)}$} In this special case $\SO(2,\C)\simeq\C^*$ and hence $\mathfrak{so}(2,\C)$ is one-dimensional and thus admits no $\mathfrak{sl}(2,\C)$-subalgebras! The maximal compact subgroup is $H=S(O(1)\times O(1))\simeq\Z_2$, so an $\SO(1,1)$-Higgs bundle can be described as a pair $(\mathcal{I},q)$ where $\mathcal{I}^2$ is a trivial line bundle and $q:\mathcal{I}\rightarrow\mathcal{I}\Ksigma$, i.e.\ $q\in H^0(\Ksigma)$. The primary topological invariant is $sw_1\in H^1(S,\Z_1)$.  It follows that the moduli space is a union of $2^{2g}$ copies of $H^0(\Ksigma)$, with the line bundle in each determined by $sw_1$.  Each component contains a point where $q=0$.  The only invariant polynomial on $\mathfrak{so}(2,\C)$ is the (degree one) Pfaffian, defined by

\begin{equation}
P\big (\begin{bmatrix}0&b\\-b&0\end{bmatrix}\big )=b\ .
\end{equation}

\noi The base of the Hitchin fibration is thus $H^0(\Ksigma)$. There are $2^{2g}$ `sections' given by the maps $q\mapsto [\mathcal{I},q]$, i.e.\ by the identification of $H^0(\Ksigma)$ with each component. These all contain points with $q=0$, i.e.\ with vanishing Higgs field.  


\subsection{Exceptional real forms}\hfill
\bigskip

The information is least complete for real forms with exceptional Lie algebras, partly because the vector bundle versions of their $G$-Higgs bundles are difficult to work with, but also because less has been worked out for their maximal compact subgroups or (in cases where they admit magical triples) for their Cayley partners.

By \cite{liegroupsteichmuller} and also \cite{BCGGO21} we know that if $G$ is a split real form of an exceptional group then $\Mhiggs G)$ contains at least one Cayley (in this case, Hitchin or Teichm\"uller) component. It also contains as many components with trivial local minima of the Hitchin function as there are topological types of principal $H$-bundles where $H$ is the maximal compact subgroup of $G$. 

For the quaternionic real forms of $F_4, E_6, E_7$ and $E_8$ the number of Cayley components is determined by the number of components in $\mathcal{M}_{\Ksigma^4}(G_{0,ss}^\R)$ in \eqref{ECayley}. This is computed in \cite{BCGGO21} for all cases except $E_8$. As stated there (see Remark 7.23) the expectation for $E_8$ is that there is just one Cayley component since the maximal compact of $G_{0,ss}^\R$ in that case (type $F_4$) is simply connected. It follows that the corresponding principal bundles have just one topological type. Moreover as stated at the beginning of this section, it is expected that in general the only connected components not labeled by topological invariants of the principal bundles are Cayley components. 

By Theorem \ref{classification} the exceptional groups with Lie algebras $\mathfrak{e}_6^2, \mathfrak{e}_7^{-5}, \mathfrak{e}_8^{-24}, \mathfrak{f}_4^4$ are those for which $\Mhiggs, G)$ have Cayley components. The number of such components depends also on the center of the group or, equivalently on their fundamental groups.  For groups of type $F_4$ or $E_8$, the adjoint forms are simply connected, so there is only one choice of group. For groups of type $E_6$ or $E_7$, the adjoint forms are not simply connected.

 \section{Tables}\label{Tables}
 \clearpage
\pagebreak


\begin{table}
\centering
  \resizebox{\textwidth}{!}{
\begin{tabular}{|c|c|c|c|c|}
\hline
&&&&\\
$G$ & $\#$ with&Cayley&Total  &References\\
& $\Phi=0$&components&$|\pi_0(\Mhiggs G))|$&\\
\hline
$\SL(2,\R)=\Sp(2,\R)=\SU(1,1)$&1&$2^{2g+1}$&$2^{2g+1}+2g-3$&\cite{goldman-88}\\
\hline
$\PSL(2,\R)=\PSp(2,\R)=\SO_0(1,2)$&1&1&$4g-3$&\cite{goldman-88}\\
\hline
$\PGL(2,\R)=\SO(1,2)$&$2^{2g+1}-1$&$(2g-2)$&$2^{2g+1}+2g-3$& \cite{so(pq)BCGGO,BS19,CollierSOnn+1components,AndrePGLnR}\\
\hline
$\SL(2m,\R),m>1$&2&$2^{2g}$&$2^{2g}+2$&\cite{liegroupsteichmuller}\\
\hline
$\PSL(2m+1,\R)$&2&1&3&\cite{liegroupsteichmuller}\\
\hline
$\PSL(2m,\R), m>1$&4&2& 6  &\cite{liegroupsteichmuller}\\
\hline
$\PGL(2m,\R), m>1$&$2^{2g+1}+1$&1&$2^{2g+1}+2$&\cite{AndrePGLnR}\\
\hline
$\SU(p,q), p\le q$&1&$2\cdot 2^{2g}$&$2\cdot 2^{2g}+2p(g-1)-1$&\cite{gothen-01,bradlow-garcia-prada-gothen:2005,BGRmaximalToledo}\\
\hline
$\PU(p,q)$&$gcd(p,q)$&$2\cdot gcd(p,q)$&$2(p+q) min\{p, q\}(g-1)$&\cite{bradlow-garcia-prada-gothen:2003,bradlow-garcia-prada-gothen:2005}\\
&&&$+ gcd(p,q)$&\cite{chains-2018,MXPupp}\\
\hline
$U^*(2m)$&1&-&1&\cite{AndreOscarSUstar}\\
\hline
\end{tabular}
}
\medskip
\caption {$\mathfrak{g}^{\C}=\mathfrak{sl}(n,\C)$}\label{SLnTable}
\end{table}

\begin{table}
\resizebox{\textwidth}{!}{
\begin{tabular}{|c|c|c|c|c|}
\hline
&&&&\\
$G$ & $\#$ with&Cayley&Total  &References\\
& $\Phi=0$&components&$|\pi_0(\Mhiggs G))|$&\\
\hline
$\Sp(4,\R)$&1&$6\cdot 2^{2g}+4g-8$&$6\cdot 2^{2g}+8g-13$&\cite{gothen-01, GM04}\\
\hline
$\PSp(4,\R)=\SO_0(2,3)$&2&$4\cdot 2^{2g}+8g-10$&$4\cdot 2^{2g}+16g-20$&\cite{AC,GO12}\\
\hline
$\Sp(2n,\R),n\ge 3$&1&$6\cdot 2^{2g}$& $[ [6\cdot 2^{2g}+2n(g-1)-1]]$&\cite{GGM}\\
\hline
$\PSp(4n+2,\R),\ n\ge 3$&1& $6$& $[[4(2n+1)(g-1)+5]]$&\cite{GP-2017}\\
 $\PSp(4n,\R),\ n\ge 4$&2&$2^{2g+2}+4$&$[[ 4(2n)(g-1)+2+2^{2g+2}]]$&\cite{GP-2017}\\
\hline
$\Sp(2p,2q)$&1&-&1&\cite{Sp(2p2q)modulispaceconnected,schaposnikPhD}\\
\hline
\end{tabular}
}
\medskip
\caption {$\mathfrak{g}^{\C}=\mathfrak{sp}(2n,\C)$}\label{Sp2nTable}
\end{table}


\begin{table}
\resizebox{\textwidth}{!}{
\begin{tabular}{|c|c|c|c|c|}
\hline
&&&&\\
$G$ & $\#$ with&Cayley&Total  &References\\
& $\Phi=0$&components&$|\pi_0(\Mhiggs G))|$&\\
\hline
$\SO(1,1)$&$2^{2g}$&0&$2^{2g}$&\cite{so(pq)BCGGO}\\
\hline
$\SO(1,2n+1),n\ge 1$&$2^{2g+1}$&-&$2^{2g+1}$&\cite{so(pq)BCGGO}\\
\hline
$\SO_0(1,2n+1)n\ge 1$&2&-&2&\cite{AG11}\\
\hline
$\SO(2,2)$&$2\cdot 2^{2g+1}-3$&$2^{2g+1}$&$6\cdot2^{2g}-3+2g(2g-3)$&\cite{so(pq)BCGGO, BSrank2}\\
\hline
$\SO(2,2n-2), n>2$&$2\cdot 2^{2g+1}-2$&$2^{2g+1}$&$[[ 6\cdot 2^{2g}+4g-8]]$&\cite{so(pq)BCGGO}\\
\hline
$\SO_0(2,2n-2), n>2$&2&$2^{2g+1}$&$[[ 2^{2g+1}+8g-10]]$&\cite{bradlow-garcia-prada-gothen:2005,BGRmaximalToledo,CTT19}\\
\hline
$\PSO_0(2,2n-2), n\ge 2$&2&$2$&$[[ 8g-6]]$&\cite{GP-2017}\\
\hline
$\SO(p,2n-p), 2<p\le n$&$4\cdot 2^{2g}$& $2\cdot 2^{2g}$&$6\cdot 2^{2g}$&\cite{so(pq)BCGGO} \\
\hline
$\SO_0(p,2n-p), 2<p\le n$&$4$& $2\cdot 2^{2g}\ (p\ \mathrm{even})$&$2\cdot2^{2g}+4\ (p\ \mathrm{even})$&\cite{so(pq)BCGGO} \\
&& $2\  (p\ \mathrm{odd})$&$6\ (p\ \mathrm{odd})$&\cite{so(pq)BCGGO} \\
\hline
$\SO^*(4n),n\ge 2$&1&2&$[[ 4n(g-1)]]$&\cite{bradlow-garcia-prada-gothen:2005,SO2n*connected}\\
$\SO^*(4n+2),n\ge1$&1&-&$[[ 4n(g-1)]]$&\cite{bradlow-garcia-prada-gothen:2005,SO2n*connected}\\
\hline
 $\PSO^*(4n), n\ge 2$&1&$2$&$[[ 2+2n(g-1)]]$&\cite{GP-2017}\\
 $\PSO^*(4n+2),n\ge 1$&1&-&$[[ 1+2(n-1)(g-1)]]$&\cite{GP-2017}\\
\hline
\end{tabular}
}
\medskip
\caption {$\mathfrak{g}^{\C}=\mathfrak{so}(2n,\C)$.}\label{SO2nTable}
\end{table}




\begin{table}
\resizebox{\textwidth}{!}{
\begin{tabular}{|c|c|c|c|c|}
\hline
&&&&\\
$G$ & $\#$ with&Cayley&Total  &References\\
& $\Phi=0$&components&$|\pi_0(\Mhiggs G))|$&\\
\hline
$\SO(1,2)=\PGL(2,\R)$&$2^{2g+1}-1$&$(2g-2)$&$2^{2g+1}+2g-3$& \cite{so(pq)BCGGO,BS19,CollierSOnn+1components,AndrePGLnR}\\
\hline
$\SO_0(1,2)=\PSL(2,\R)$&1&1&$4g-3$&\cite{goldman-88}\\
\hline
$\SO(1,2n), n>1$&$2^{2g+1}$&-&$2^{2g+1}$&\cite{so(pq)BCGGO}\\
\hline
$\SO_0(1,2n)n\ge 1$&2&-&2&\cite{AG11}\\
\hline
$\SO(2,3)$&$4\cdot 2^{2g}-2$&$2\cdot 2^{2g}+4g-5$&$6\cdot 2^{2g}+8g-13$&\cite{GO12,CollierSOnn+1components}\\
\hline
$\SO_0(2,3)=P\mathrm{Sp}(4,\R)$&2&$4\cdot 2^{2g}+8g-10$&$4\cdot 2^{2g}+16g-20$&\cite{AC,GO12}\\
\hline
$\SO(2,2n-1),n>2$&$2\cdot 4^{2g}-2$&$2\cdot 2^{2g}$&$[[ 6\cdot 2^{2g}+4g-8]]$&\cite{so(pq)BCGGO}\\
\hline
$\SO_0(2,2n-1)=\PSO_0(2, 2n-1)$&2&$2^{2g+1}$&$[[ 2^{2g+1}+8g-10]]$&\cite{bradlow-garcia-prada-gothen:2005,BGRmaximalToledo}\\
$n>2$&&&&\cite{CTT19}\\
\hline
$\SO(p,2n+1-p),2<p<n$&$4\cdot 2^{2g}$& $2\cdot 2^{2g}$&$6\cdot 2^{2g}$&\cite{so(pq)BCGGO} \\
\hline
$\SO_0(p,2n+1-p), 2<p< n$&$4$& $2\cdot 2^{2g}\ (p\ \mathrm{even})$&$2\cdot2^{2g}+4\ (p\ \mathrm{even})$&\cite{so(pq)BCGGO} \\
&& $2\  (p\ \mathrm{odd})$&$6\ (p\ \mathrm{odd})$&\cite{so(pq)BCGGO} \\
\hline
$\SO(n,n+1),n>2$&$4\cdot 2^{2g}$&$2\cdot 2^{2g}$&$6\cdot 2^{2g}$&\cite{CollierSOnn+1components}\\
&&$+2n(g-1)-1$&$+2n(g-1)-1$&\\
\hline
\end{tabular}
}
\medskip
\caption {$\mathfrak{g}^{\C}=\mathfrak{so}(2n+1,\C)$}\label{SO2n+1Table}

\end{table}

 
\begin{table}
\resizebox{\textwidth}{!}{
\begin{tabular}{|c|c|c|c|c|}
\hline
&&&&\\
$G$ & $\#$ with&Cayley&Total  &References\\
& $\Phi=0$&components&$|\pi_0(\Mhiggs G))|$&\\
\hline
$E^6_6$ &&[[1]]&&\cite{liegroupsteichmuller,BCGGO21}\\
\hline 
$E^2_6$ &&1 (simply connected)&&\cite{BCGGO21}\\
&&3 (adjoint))&&\cite{BCGGO21}\\
\hline 
$E^{-14}_6$ &&2&&\cite{BGRmaximalToledo,GP-2017}\\
\hline 
$E^{-26}_6$ &[[1]]&-&[[1]]&\cite{collier2023}\\
\hline
$E^7_7$&&[[1]]&&\cite{liegroupsteichmuller,BCGGO21}\\
\hline
$E_7^{-25}$&&$>1$&&\cite{BGRmaximalToledo}\\
\hline 
$E^{-5}_7$&&1 (simply connected)&&\cite{BCGGO21}\\
&&2 (adjoint))&&\cite{BCGGO21}\\
\hline
$E_8^8$&&[[1]]&&\cite{liegroupsteichmuller,BCGGO21}\\
\hline
$E^{-24}_8$& &[[1]]&&\cite{BCGGO21}\\
\hline
$F^{-20}_4$ &[[1]]&-&[[1]]&\cite{collier2023}\\
\hline
$F^4_4$ &1&3&4&\cite{liegroupsteichmuller,BCGGO21}\\
\hline
$G^2_2$ &&[[1]]&&\cite{liegroupsteichmuller,BCGGO21}\\
\hline
\end{tabular}
}
\medskip
\caption {Exceptional $\mathfrak{g}^{\C}$}\label{Lie-excep}
\end{table}

\clearpage


\bibliographystyle{amsplain}

\begin{thebibliography}{99}

\bibitem{allesandrini2019} Daniele Alessandrini, {\it Higgs Bundles and Geometric Structures on Manifolds}, SIGMA 15 (2019)

\bibitem{AMTW} Daniele Alessandrini, Sara Maloni, Nicolas Tholozan, Anna Wienhard, {\it Fiber bundles associated with Anosov representations}, arXiv:2303.10786 (2023)

\bibitem{AC}Daniele Alessandrini and Brian Collier {\it The geometry of maximal components of the $P\mathrm{Sp}(4,\R)$- character variety} Geometry \& Topology 23 (2019) 1251-1337

\bibitem{AG11}Marta Aparicio Arroyo and Oscar Garcia-Prada, \newblock Higgs bundles for the Lorentz group, \newblock {\it Illinois Journal of Mathematics} 55 (2011) 1299--1326.

\bibitem{so(pq)BCGGO}
Marta Aparicio-Arroyo, Steven Bradlow, Brian Collier, Oscar Garc\'{\i}a-Prada,
  Peter~B. Gothen, and Andr\'{e} Oliveira.
\newblock {$\mathrm{SO}(p,q)$}-{H}iggs bundles and higher {T}eichm\"{u}ller
  components.
\newblock {\em Invent. Math.}, 218(1):197--299, 2019.

\bibitem{BS19} David Baraglia and Laura Schaposnik, {\it Cayley and Langlands-type correspondences for orthogonal Higgs bundles} Trans. Am. Math. Soc. 371(10), 7451--7492 (2019)

\bibitem{BSrank2}
David Baraglia and Laura Schaposnik.
\newblock Monodromy of rank 2 twisted {H}itchin systems and real character
  varieties.
\newblock {\em Trans. Amer. Math. Soc.}, 370(8):5491--5534, 2018.

\bibitem{DavidLauraCayleyLanglands}
David Baraglia and Laura~P. Schaposnik.
\newblock Cayley and {L}anglands type correspondences for orthogonal {H}iggs
  bundles.
\newblock {\em Trans. Amer. Math. Soc.}, 371(10):7451--7492, 2019.

\bibitem{BGRmaximalToledo}
Olivier Biquard, Oscar Garcia-Prada, and Roberto Rubio.
\newblock Higgs bundles, the {T}oledo invariant and the {C}ayley
  correspondence.
\newblock {\em Journal of Topology}, 10(3):795--826, 2017.

\bibitem{BCGGO21}
 Steven Bradlow, Brian Collier, Oscar Garc\'{\i}a-Prada,
  Peter~B. Gothen, and Andr\'{e} Oliveira.
{\it A general Cayley correspondence and higher Teichm\"uller  spaces}, arXiv:2101.09377 (2021)
\url{https://urldefense.com/v3/__https://arxiv.org/pdf/2101.09377.pdf__;!!DZ3fjg!6ba3Hh-fH55lNShG2t7rnZSBI3Pj3ycb9TKA5DPE7XP3mtSNtqqYcgOOd8R73VvQdV7de_umqfSfkLI0xodoGg$ }

\bibitem{bradlow-garcia-prada-gothen:2003}
S.~B. Bradlow, O.~Garc{\'\i}a-Prada,  and P.~B. Gothen, Surface group
representations and $\mathrm{U}(p,q)$-{H}iggs bundles,
{\em J. Differential Geom.}
\textbf{64} (2003), 111--170.

\bibitem{bradlow-garcia-prada-gothen:2004}
S.~B. Bradlow, O.~Garc{\'\i}a-Prada,  and P.~B. Gothen,  Representations of surface groups in the general linear group, Publicaciones de la RSME, Vol. 7 (2004) 83--94.


\bibitem{bradlow-garcia-prada-gothen:2005}
S.~B. Bradlow, O.~Garc{\'\i}a-Prada,  and P.~B. Gothen, Maximal surface group
representations in isometry groups of classical Hermitian symmetric spaces,
{\em Geometria Dedicata.}
\textbf{122} (2006), 185--213.

\bibitem{bradlow-garcia-prada-gothen:2008}
S.~B. Bradlow, O.~Garc{\'\i}a-Prada,  and P.~B. Gothen, {\it Homotopy groups of moduli spaces of representations} Topology, Volume 47, Issue 4, September 2008, Pages 203--224


\bibitem{bradlow-garcia-prada-gothen:2012}
S.~B. Bradlow, O.~Garc{\'\i}a-Prada,  and P.~B. Gothen, {\it Deformations of maximal representations in $\mathrm{Sp}(4,\R)$}  Quart. J. Math. 63 (2012), 795--843

\bibitem{chains-2018}
Steven Bradlow, Oscar Garcia-Prada, Peter Gothen, and Jochen Heinloth.
\newblock Irreducibility of moduli of semistable chains and applications to
  {$\mathrm{U}(p,q)$}-{H}iggs bundles.
\newblock In {\em Geometry and Physics: Volume 2, A Festschrift in honour of
  Nigel Hitchin}, pages 455--470. Oxford University Press, 2018.

\bibitem{HermitianTypeHiggsBGG}
Steven~B. Bradlow, Oscar Garc{\'{\i}}a-Prada, and Peter~B. Gothen.
\newblock Maximal surface group representations in isometry groups of classical
  {H}ermitian symmetric spaces.
\newblock {\em Geom. Dedicata}, 122:185--213, 2006.

\bibitem{SO2n*connected}
Steven~B. Bradlow, Oscar Garc{\'{\i}}a-Prada, and Peter~B. Gothen.
\newblock Higgs bundles for the non-compact dual of the special orthogonal
  group.
\newblock {\em Geom. Dedicata}, 175:1--48, 2015.


\bibitem{collier2023} Brian Collier, {\it Classification of character varieties with no higher Teichm\"uller components} (2023) (in preparation)

\bibitem{CollierSOnn+1components}
Brian Collier.
\newblock {$\mathsf{SO}(n,n+1)$-surface group representations and their Higgs
  bundles}.  Annales scientifiques de l'ENS volume 53, issue 6 (2020) 1561--1616.

\bibitem{CollierHitchin} Brian Collier
\newblock {Various generalizations and deformations of $\PSL(2,\R)$ surface group representations and their Higgs bundles}. Geometry and Physics: a Festschrift in honour of Nigel Hitchin, Oxford Press.

\bibitem{CTT19} B. Collier,  N. Tholozan and J. Toulisse, {\it The geometry of maximal representations of surface groups into $SO(2,n)$},  Duke Math Journal, Volume 168, Number 15 (2019), 2873-2949.


\bibitem{garcia-pradasurvey09}O. Garcia-Prada, {\it Higgs bundles and surface group representations, in Moduli Spaces and Vector Bundles}, LMS Lecture Notes Series 359, 265--310, Cambridge University Press, [2009].

\bibitem{garcia-pradasurvey}O. Garcia-Prada {\it Higgs bundles and higher Teichm\"uller  spaces}, in Handbook on Teichm\"uller  theory Athanase Papadopoulos (editor), Vol. VII, 239?285, IRMA Lect. Math. Theor. Phys., 30, Eur. Math. Soc., Zurich,[2020]


\bibitem{GGM}
Oscar Garc{\'{\i}}a-Prada, Peter Gothen, and Ignasi Mundet~i Riera.
\newblock Higgs bundles and surface group representations in the real symplectic
  group.
\newblock {\em Journal of Topology}, 6(1):64--118, 2013.

\bibitem{Oliveira_GarciaPrada_2016}
Oscar Garc{\'{\i}}a-Prada and Andr{\'e} Oliveira.
\newblock Connectedness of {H}iggs bundle moduli for complex reductive {L}ie
  groups.
\newblock {\em Asian Journal of Mathematics}, 21(5):791--810, 2017.

\bibitem{AndreOscarSUstar}
Oscar Garc\'{\i}a-Prada and Andr\'{e}~G. Oliveira.
\newblock Higgs bundles for the non-compact dual of the unitary group.
\newblock {\em Illinois J. Math.}, 55(3):1155--1181 (2013), 2011.

\bibitem{Sp(2p2q)modulispaceconnected}
Oscar Garc{\'{\i}}a-Prada and Andr{\'e}~G. Oliveira.
\newblock Connectedness of the moduli of {${\rm Sp}(2p,2q)$}-{H}iggs bundles.
\newblock {\em Q. J. Math.}, 65(3):931--956, 2014.

\bibitem{GP-2017}
Oscar Garcia-Prada, Andre Oliveira, {\it Maximal Higgs bundles for adjoint forms via Cayley correspondence} Geometriae Dedicata, 190, No.1 (2017), 1--22

\bibitem{GM04}O. Garcia-Prada, I. Mundet i Riera, {\it Representations of the fundamental group of a closed oriented surface in $Sp(4, \R)$}, Topology 43 (2004), 831--855.


\bibitem{goldman-88}
William~M. Goldman.
\newblock Topological components of spaces of representations.
\newblock {\em Invent. Math.}, 93(3):557--607, 1988.

\bibitem{gothen-01}
Peter~B. Gothen.
\newblock Components of spaces of representations and stable triples.
\newblock {\em Topology}, 40(4):823--850, 2001.

\bibitem{gothensurvey}P.B. Gothen, {\it Representations of surface groups and Higgs bundles}, in Moduli Spaces, Edited by L. Brambila-Paz, O. Garcia-Prada, P. Newstead and R.P.  Thomas, LMS Lecture Note Series, 411, CUP, 2014. 


\bibitem{GO12}Peter B. Gothen and Andre G. Oliveira. {\it Rank two quadratic pairs and surface group representations} Geom. Dedicata, 161:335--375, 2012.

\bibitem{guichardsurvey}O. Guichard, {\it An Introduction to the Differential Geometry of Flat Bundles and of Higgs Bundles}, Lecture Notes Series, Institute for Mathematical Sciences, National University of Singapore, The Geometry, Topology and Physics of Moduli Spaces of Higgs Bundles, pp. 1--63 (2018)


\bibitem{GW10} O. Guichard and A. Wienhard, {\it Topological invariants of Anosov representations}, Journal of Topology 3 (2010), 578--642.

\bibitem{GW18} O. Guichard and A. Wienhard, {\it  Positivity and higher Teichm\"uller theory.} In: Proceedings of
the 7th European Congress of Mathematics, pp. 289--310. European Mathematical Society,
Zurich (2018)

\bibitem{GW22} O. Guichard and A. Wienhard, {\it Generalizing Lusztig's total positivity}, arXiv:2208.10114  (2022)

\bibitem{Helgason} S. Helgason, {\it Differential geometry, Lie groups, and symmetric spaces}, volume 34 of Graduate Studies in Mathematics. American Mathematical Society, Providence, RI, (2001). [Corrected reprint of the 1978 original.]

\bibitem{selfduality}
Nigel Hitchin.
\newblock The self-duality equations on a {R}iemann surface.
\newblock {\em Proc. London Math. Soc. (3)}, 55(1):59--126, 1987.

\bibitem{IntSystemFibration}
Nigel Hitchin.
\newblock Stable bundles and integrable systems.
\newblock {\em Duke Math. J.}, 54(1):91--114, 1987.

\bibitem{liegroupsteichmuller}
Nigel Hitchin.
\newblock Lie groups and {T}eichm\"uller space.
\newblock {\em Topology}, 31(3):449--473, 1992.

\bibitem{Kn}A. Knapp, {\it Lie groups beyond an introduction}, Birkhauser, Progress in Mathematics vol. 140, (1996)

\bibitem{K}B. Kostant {\i Lie group representations on polynomial rings}, Amer. J. Math., 85, 327--404 (1963)

\bibitem{Li2019}Qi. Li  {\it An Introduction to Higgs Bundles via Harmonic Maps}, SIGMA 15 (2019) [Special Issue on Geometry and Physics of Hitchin Systems ]


\bibitem{JunLiConnectedness}
Jun Li.
\newblock The space of surface group representations.
\newblock {\em Manuscripta Math.}, 78(3):223--243, 1993.

\bibitem{MXPupp}E. Markman and E. Z. Xia. {\it The moduli of flat $\PU(p,p)$-structures with large Toledo invariants} Math. Z., 240(1):95--109, 2002.

\bibitem{AndrePGLnR}
Andr\'e~Gama Oliveira.
\newblock Representations of surface groups in the projective general linear
  group.
\newblock {\em Internat. J. Math.}, 22(2):223--279, 2011.

Comments: v3  (\url{https://urldefense.com/v3/__https://arxiv.org/pdf/0901.2314v3.pdf__;!!DZ3fjg!6ba3Hh-fH55lNShG2t7rnZSBI3Pj3ycb9TKA5DPE7XP3mtSNtqqYcgOOd8R73VvQdV7de_umqfSfkLK3JyF09w$ } ) included an erratum (Section 12) which shows why Theorem 1.3 (stating that the Hitchin component in PSL(3,R) is homotopically equivalent to PSO(3)) is not correct, even though the original manuscript is left unchanged. This erratum has been published in Int. J. Math., 30, No. 2 (2019) 

\bibitem{schaposnikPhD}L. P. Schaposnik, {\it Spectral data for G-Higgs bundles}, D. Phil. Thesis, University of Oxford, 2012.

\bibitem{schaposniksurvey}L. Schaposnik,  {\it Higgs bundles - recent applications}. Notices Amer. Math. Soc. 67 (2020), no. 5, 625--634.


\bibitem{schmitt}A. Schmitt, {\it Geometric Invariant Theory and Decorated Principal Bundles} Zurich Lectures in Advanced Mathematics, EMS (2000)

\bibitem{Sim88}C. T. Simpson, {\it Moduli of representations of the fundamental group of a smooth projective variety. I}, Publ. Math. Inst. Hautes Etudes Sci. 79 (1994), 47--129.

\bibitem{swobodasurvey}J. Swoboda, {\it Moduli Spaces of Higgs Bundles ? Old and New} , Jahresber Dtsch Math-Ver  Jahresber. Dtsch. Math.-Ver. 123 (2021), no. 2, 65--130. 

\bibitem{wentworthsurvey}R. Wentworth. {\it Higgs bundles and local systems on Riemann surfaces} 
In Geometry and Quantization of Moduli Spaces, CRM Advanced Courses in Mathematics, Birkhauser/Springer, 2016.

\bibitem{W2018}Anna Wienhard. {\it An invitation to higher Teichm\"uller  theory.} In Proceedings of the International Congress of Mathematicians, Rio de Janeiro 2018. Vol. II. Invited lectures, pages 1013--1039. World Sci. Publ., Hackensack, NJ, (2018)

\bibitem{Xia} E. Z. Xia, {\it The moduli of flat $PGL(2,\R)$ connections on Riemann surfaces}, Commun. Math. Phys. 203 (1999), 531--549.
\end{thebibliography}

\end{document}